\documentclass[12pt,reqno]{amsart}
\usepackage{amssymb}
\usepackage{graphicx}

\usepackage{amsmath}
\usepackage{a4wide}
\usepackage{epsf}
\usepackage{epsfig}
\usepackage{psfrag}

\newtheorem{lemma}{Lemma}
\newtheorem{theorem}{Theorem}
\newtheorem{corollary}[lemma]{Corollary}
\newtheorem{proposition}[lemma]{Proposition}
\def\om{\omega}

\title[Cyclotomic matrices]{Integer symmetric matrices having all their eigenvalues in the interval $[-2,2]$}
\author[]{James McKee}
\address{Department of Mathematics\\
Royal Holloway, University of London\\
Egham Hill\\
Egham\\
Surrey TW20 0EX\\
UK}
\author[]{Chris Smyth}
\address{School of Mathematics and Maxwell Institute for Mathematical Sciences \\
University of Edinburgh\\
James Clerk Maxwell Building\\
King's Buildings\\
Mayfield Road\\
Edinburgh  EH9 3JZ\\
UK}

\begin{document}
\begin{abstract}
We completely describe all integer symmetric matrices that have
all their eigenvalues in the interval $[-2,2]$. Along the way we
classify all signed graphs, and then all charged signed graphs,
having all their eigenvalues in this same interval. We then
classify subsets of the above for which the integer symmetric
matrices, signed graphs and charged signed graphs have all their
eigenvalues in the open interval $(-2,2)$.
\end{abstract}
\maketitle

\section{Introduction}
Let $A$ be an $n\times n$ integer symmetric matrix with
characteristic polynomial $\chi_A(x)=\det(xI-A)$. The aim of this
paper is to describe all such matrices $A$ that have the maximum
modulus of their eigenvalues at most $2$.  The significance of the
bound $2$ is that, by a result of Kronecker \cite{K},  every
eigenvalue of such a matrix $A$ is then of the form $\om+\om^{-1}$, for
some root of unity $\om$. Thus $z^n\chi_A(z+1/z)$ is a cyclotomic
polynomial. For this reason we call such integer symmetric
matrices {\it cyclotomic} matrices.

In 1970 J.H.~Smith \cite{Smi} classified all cyclotomic
$\{0,1\}$-matrices   with zeros on the diagonal, regarding them as
adjacency matrices of graphs (see Figure \ref{F-except}). Such graphs were called
\emph{cyclotomic} graphs in \cite{MS}. It turns out that a full
description of cyclotomic matrices is conveniently stated using
more general graphs. So if we allow the off-diagonal elements of
our matrix to be chosen from the set $\{-1,0,1\}$, we obtain a
{\it signed} graph (see \cite{CST},\cite{Z2}),  a non-zero
$(i,j)$th entry  denoting a `sign' of $-1$ or $1$ on the edge
between vertices $i$ and $j$. Further, for a general symmetric
$\{-1,0,1\}$ matrix, where now the diagonal entries may be
nonzero, we obtain what we call a {\it charged} signed graph; we
regard a nonzero $(i,i)$th entry of $A$ as corresponding to a
`charge' on its $i$th vertex. If none of the edges of a charged
signed graph in fact have sign $-1$, then we have a charged
(unsigned) graph. However, a graph is also a signed graph, and a
signed graph is also a charged signed graph. The notion of a
charged signed graph is a convenient device for picturing and
discussing symmetric integer matrices with entries in
$\{-1,0,1\}$. These are the most important matrices in our
description of general cyclotomic matrices.

In this paper we extend Smith's result to cyclotomic charged
signed graphs (Theorem \ref{T-LCSG}), and then, with little
further work, to all cyclotomic matrices (Theorem \ref{T-ISM}).
Along the way we find all cyclotomic signed graphs (Theorem
\ref{T-CSG}). As a consequence, we can also describe all
cyclotomic charged graphs (Theorem   \ref{T:unsigned_charged}) and
all cyclotomic matrices whose entries are non-negative (Theorem
\ref{T:nonnegative_matrix}).

Having obtained our results for the closed interval $[-2,2]$, it
is then very natural to consider restricting the eigenvalues to
the open interval $(-2,2)$. We give a complete classification of
symmetric integer matrices with eigenvalues in this restricted set
(Theorem \ref{T-ISMopen}). As in the case of the closed interval,
there are corresponding results for cyclotomic signed graphs
(Theorem \ref{T:open_signed}),  cyclotomic charged signed graphs
(Theorem \ref{T:open_charged_signed}), cyclotomic charged graphs
(Theorem  \ref{T:open_charged}) and cyclotomic matrices whose
entries are non-negative (Theorem
\ref{T:nonnegative_matrix_open}). Having dealt with the general
cyclotomic case, this is a relatively straightforward problem.
There is a connection here with the theory of finite reflection
groups and their Coxeter graphs, and we conclude with a discussion
of this.

In \cite{MS}, cyclotomic graphs were used to construct Salem
numbers and Pisot numbers. The original motivation for this
current work was that it provides one of the ingredients necessary
to extend the work in \cite{MS}. But we think that our results may
be of independent interest.

Throughout the paper, a {\it subgraph} of the (charged, signed)
graph under consideration will always mean a vertex-deleted
subgraph, that is, an induced subgraph on a subset of the
vertices.

\section{Interlacing, and reduction to maximal indecomposable matrices}
In order to state our results, we need some preliminaries. The
matrix $A$ will be called \emph{indecomposable} if and only if the
underlying  graph is \emph{connected}. (In the underlying graph,
vertices $i$ and $j$ are adjacent if and only if the $(i,j)$th
entry of $A$ is nonzero.) If $A$ is not indecomposable, then there
is a reordering of the rows (and columns) such that the matrix has
block diagonal form with more than one block, and its list of
eigenvalues is found by pooling the lists of the eigenvalues of
the blocks. For our classification of cyclotomic  matrices, it is
clearly sufficient to consider indecomposable ones.

A repeatedly useful tool for us is Cauchy's interlacing theorem
(for a short proof, see \cite{Fis}).

\begin{lemma}[Interlacing Theorem]\label{L:interlacing}
Let $A$ be a real symmetric matrix, with eigenvalues
$\lambda_1 \le \lambda_2 \le \ldots \le \lambda_n$.
Pick any row $i$, and let $B$ be the matrix formed
by deleting row $i$ and column $i$ from $A$.
Then the eigenvalues of $B$ interlace with those of $A$:
if $B$ has eigenvalues $\mu_1 \le \ldots \le \mu_{n-1}$,
then
\[
\lambda_1 \le \mu_1 \le \lambda_2 \le \mu_2 \le \ldots \le \mu_{n-1} \le \lambda_n\,.
\]
\end{lemma}

In view of this Lemma, if $A$ is cyclotomic, then so
is any matrix obtained by deleting from $A$ any number
of its rows, along with the corresponding columns:
we then speak of the smaller matrix as being \emph{contained in}
the larger one (the smaller graph is an induced subgraph of the larger graph).
We call an indecomposable cyclotomic matrix (or its graph)
\emph{maximal} if it
is not contained in a strictly-larger indecomposable cyclotomic matrix:
the corresponding cyclotomic graph is not an induced
subgraph of a strictly larger connected cyclotomic graph.
We shall see that every non-maximal indecomposable cyclotomic
matrix is contained in a maximal one.
It is therefore enough for us to classify all
maximal indecomposable cyclotomic matrices.

When we consider matrices that have all their eigenvalues in the
open interval $(-2,2)$, we shall see that it is no longer always
true that every such matrix is contained in a maximal one: there
is an infinite family of indecomposable exceptions.

\section{Equivalence, strong equivalence and switching}

 Denote by $O_n(\mathbb Z)$ the orthogonal group of $n\times n$ signed permutation matrices.
Then conjugation of a cyclotomic matrix by a matrix in
$O_n(\mathbb Z)$ gives a cyclotomic matrix with the same
eigenvalues. We say that two  $n\times n$  cyclotomic matrices are
{\it strongly equivalent} if they are related in this way.
Further, we say that two indecomposable cyclotomic matrices $A$
and $A'$ are merely {\it equivalent} if $A'$ is strongly
equivalent to $A$ or $-A$. This notion then extends easily to
decomposable cyclotomic matrices. Both of these notions are
equivalence relations on the set of all cyclotomic matrices. For
indecomposable cyclotomic matrices, the equivalence classes for
the weaker notion are the union of one or two strong equivalence
classes, depending on whether or not $-A$ is in the same strong
equivalence class as $A$.  It is clearly sufficient to classify
all cyclotomic  matrices up to equivalence.

For a charged signed graph, the notions of strong equivalence and
equivalence of course carry over via the adjacency matrix. Now
$O_n(\mathbb Z)$ is generated  by diagonal matrices of the form
$\text{diag}(1,1,\dots,1,-1,1,\dots,1)$ and by permutation
matrices.  Conjugation by these diagonal matrices corresponds to
reversing the signs of all edges incident at a certain vertex $v$;
we call this {\it  switching at $v$}. Conjugation by a permutation
matrix merely means that we can ignore vertex labels; we therefore
do not label the vertices of our graphs. Thus for unlabelled
charged signed graphs, strong equivalence classes are generated
only by such switching operations. The concept of switching, and
signed switching classes, appeared earlier for signed graphs in
\cite{CST}.

Equivalence of charged signed graphs is generated both by
switching, and by the operation of reversing all the edge signs
and  vertex charges of a component of a graph.

Since most of our graphs will in fact be signed graphs, we avoid
clutter by drawing edges with sign $1$ as unbroken lines
---------, and edges with sign $-1$ as dashed lines - - - - - -.
For vertices, those of charge $1,0,-1$ will be drawn
\textcircled{+},$\bullet$,\textcircled{--} respectively, with the
vertices $\bullet$ without a charge
  being called {\it neutral} vertices.

\section{Main results}

\begin{figure}[h]
\begin{center}
\leavevmode
\psfragscanon
\psfrag{A}{$A$}
\psfrag{B}{$B$}
\psfrag{a}{$A$}
\psfrag{b}{$B$}
\hbox{ \epsfxsize=3.3in
\epsffile{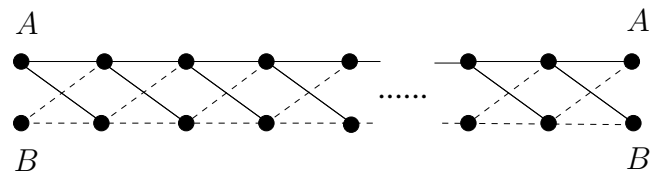} }
\end{center}
\caption{The family $T_{2k}$ of $2k$-vertex maximal connected
cyclotomic toral tesselations, for $k\ge 3$. (The two copies of
vertices A and B should be identified, as in Figure \ref{F-torus}
below.)} \label{F-CSG_general}
\end{figure}

\begin{figure}[h]
\begin{center}
\leavevmode
\hbox{
\epsfxsize=2.5in
\epsffile{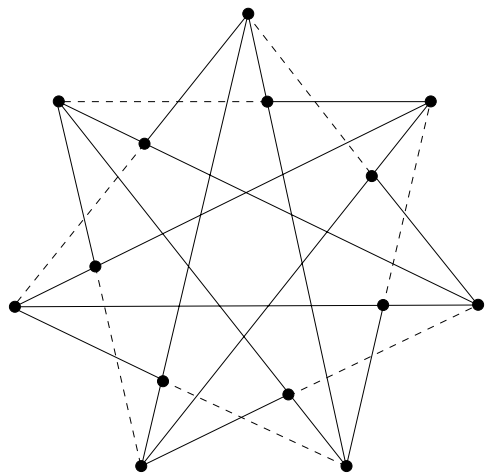}
}
\end{center}
\caption{The $14$-vertex sporadic maximal connected cyclotomic
signed graph $S_{14}$. See also Section \ref{S-Chapman}.}
\label{F-CSG_14}
\end{figure}

\begin{figure}[h]
\begin{center}
\leavevmode
\hbox{
\epsfxsize=3.3in
\epsffile{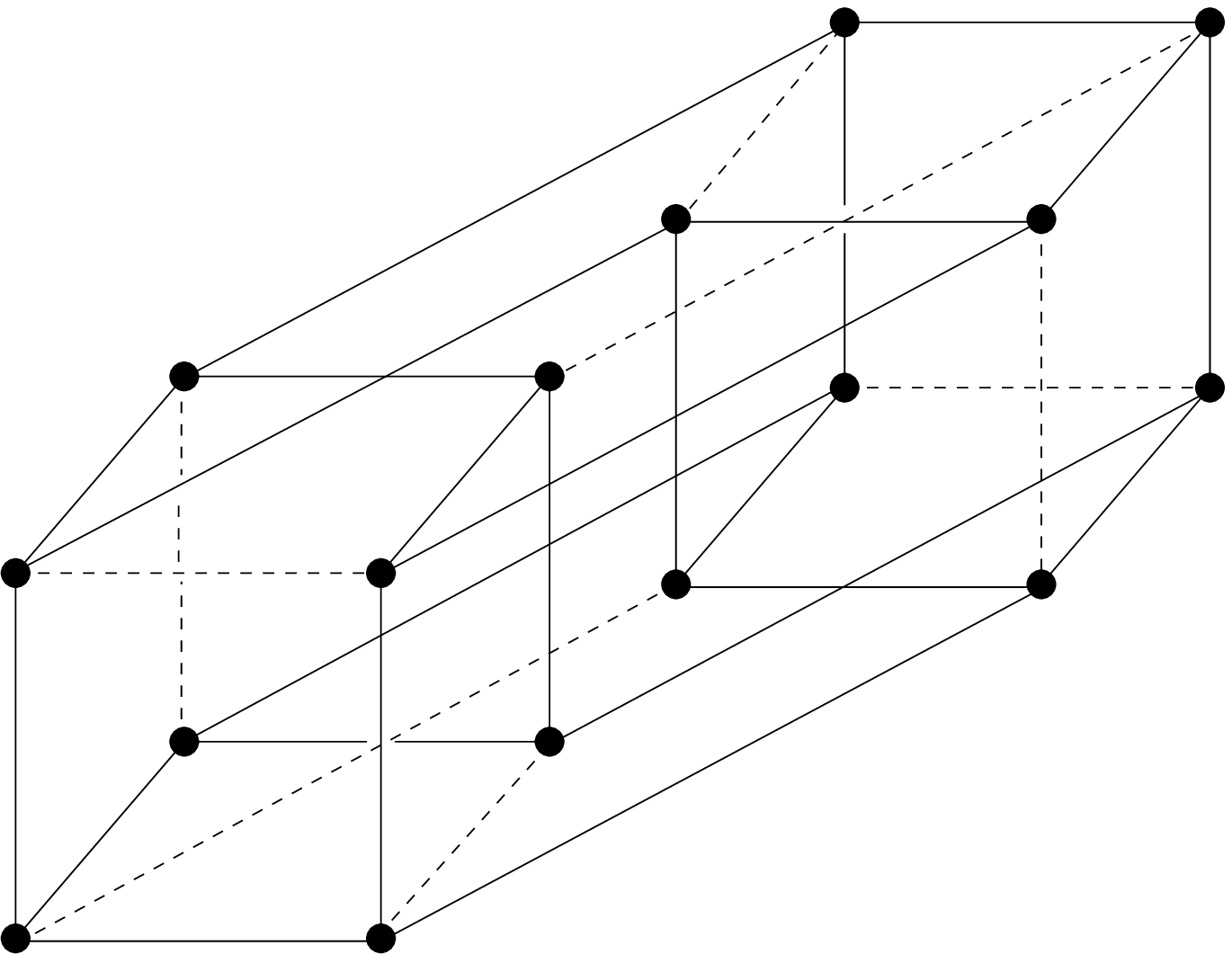}
}
\end{center}
\caption{The hypercube sporadic maximal connected cyclotomic
signed graph $S_{16}$.} \label{F-CSG_16}
\end{figure}

\begin{theorem}[{``Uncharged, signed, $[-2,2]$''}]\label{T-CSG}
Every maximal connected cyclotomic signed graph is equivalent to
one of the following:
\begin{enumerate}
\item[(i)]  For some $k=3,4,\dots$, the $2k$-vertex toral
tesselation $T_{2k}$ shown in Figure \ref{F-CSG_general};
\item[(ii)] The $14$-vertex signed graph  $S_{14}$ shown in Figure \ref{F-CSG_14};
\item[(iii)] The $16$-vertex signed hypercube $S_{16}$ shown in Figure \ref{F-CSG_16}.
\end{enumerate}
Further, every connected cyclotomic signed graph is contained in a
maximal one.
\end{theorem}

\begin{figure}[h]
\begin{center}
\epsfig{file=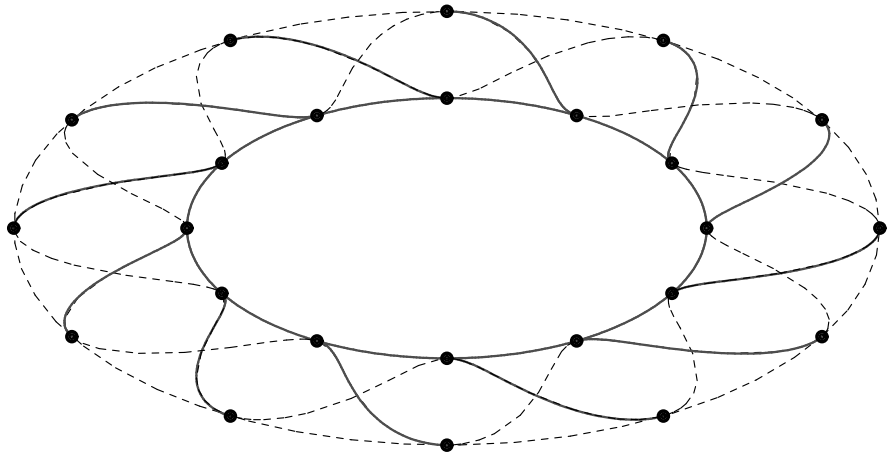,height=6cm,width=12cm}
\end{center}
\caption{A typical toral tesselation $T_{2k}$: the signed graph $T_{24}$.} \label{F-torus}
\end{figure}


\begin{figure}[h]
\begin{center}
\leavevmode
\hbox{
\epsfxsize=1.7in
\epsffile{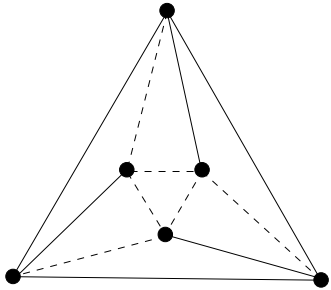}
}
\end{center}
\caption{The octahedral maximal connected cyclotomic signed graph
$T_6$.} \label{F-oct}
\end{figure}

In particular, $k=3$ of case (i) gives an octahedron $T_6$, shown in
Figure \ref{F-oct}, while a more typical example $T_{24}$ is shown in Figure \ref{F-torus}.

\begin{figure}[h]
\begin{center}
\leavevmode
\psfragscanon
\psfrag{+}{$C_{2k}^{++}$}
\psfrag{-}{$C_{2k}^{+-}$}
\hbox{ \epsfxsize=4.3in
\epsffile{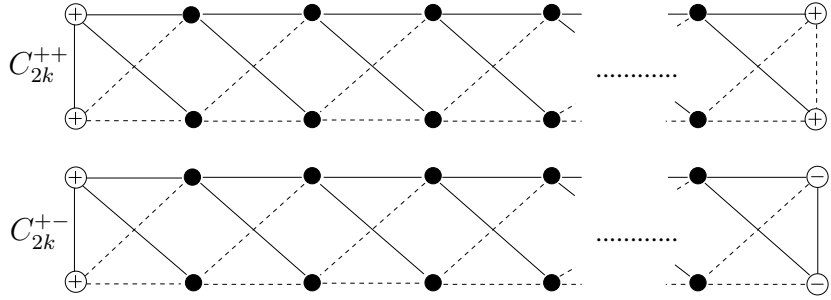} }
\end{center}
\caption{The families of $2k$-vertex maximal connected cyclotomic
cylindrical tesselations $C_{2k}^{++}$ and $C_{2k}^{+-}$, for
$k\ge 2$.} \label{F-LCSG_general}
\end{figure}

\begin{figure}[h]
\begin{center}
\leavevmode
\psfragscanon
\psfrag{7}{$S_7$}
\psfrag{8}{$S_8$}
\psfrag{9}{$S'_8$}
\hbox{ \epsfxsize=5.3in
\epsffile{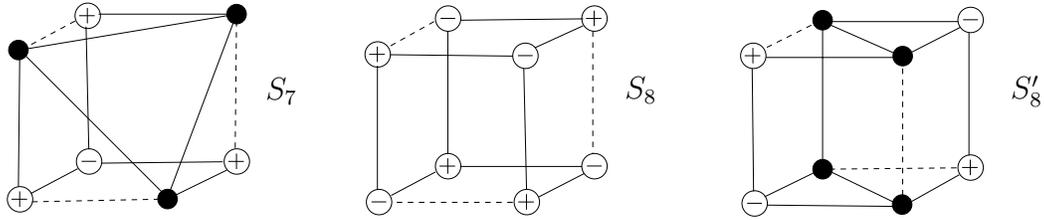} }
\end{center}
\caption{The three sporadic maximal connected cyclotomic charged
signed graphs $S_7,S_8,S'_8$.} \label{F-LCSG_sporadic}
\end{figure}

\begin{theorem}[{``Charged, signed, $[-2,2]$''}]\label{T-LCSG}
Every maximal connected cyclotomic charged signed graph not
included in Theorem \ref{T-CSG} is equivalent to one of the
following:
\begin{enumerate}
\item[(i)]  For some $k=2,3,4,\dots$, one of the two $2k$-vertex
cylindrical tesselations $C_{2k}^{++},C_{2k}^{+-}$ shown in Figure
\ref{F-LCSG_general}; \item[(ii)] One of the three sporadic
charged signed graphs $S_7,S_8,S'_8$ shown in Figure
\ref{F-LCSG_sporadic};

\end{enumerate}
Further, every connected cyclotomic charged signed graph is
contained in a maximal one.
\end{theorem}

\begin{figure}[h]
\begin{center}
\leavevmode
\psfragscanon
\psfrag{++}{$C_4^{++}$}\psfrag{+-}{$C_4^{+-}$}
\hbox{
\epsfxsize=2.3in
\epsffile{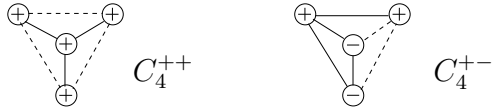}
}
\end{center}
\caption{The two maximal connected cyclotomic charged signed
tetrahedra $C_4^{++}$ and $C_4^{+-}$.} \label{F-LCSG_tet}
\end{figure}

In particular, $k=2$ of case (i) gives two charged tetrahedra
$C_4^{++},C_4^{+-}$, shown in Figure \ref{F-LCSG_tet}.

We remark that all the maximal cyclotomic graphs of Theorems
\ref{T-CSG} and \ref{T-LCSG} are `visibly' cyclotomic: their
adjacency matrices $A$ all satisfy $A^2=4I$, so all their
eigenvalues are $\pm 2$. The exact multiplicity of these eigenvalues
is given in Table \ref{Ta-1} at the end of the paper.

Our most general result is readily deduced from the previous two theorems.

\begin{theorem}[{``Integer matrix, $[-2,2]$''}]\label{T-ISM} Every maximal indecomposable  cyclotomic
matrix is equivalent to one of the following:
\begin{enumerate}
\item[(i)]  The adjacency matrix of a maximal connected  charged
cyclotomic signed graph (given by Theorems \ref{T-CSG} and
\ref{T-LCSG});

\item[(ii)]  The $1\times 1$ matrix $(2)$ or the matrix
$\left(\begin{matrix}0 & 2\\2 & 0\end{matrix}\right)$.

Further, every indecomposable cyclotomic matrix
is contained in a maximal one.
\end{enumerate}
\end{theorem}

\section{Simplifications}

A  signed graph $G$ is called \emph{bipartite} if its vertices can
be split into two disjoint parts   such that every edge of $G$
joins a vertex in one part to a vertex in the other (\cite{Z3}).
The eigenvalues of $G$ are then symmetric about 0, counted with
multiplicity; we record this fact as a Lemma.
\begin{lemma}\label{L:bipartite}
Let $G$ be a bipartite signed graph with $n$ vertices. Then
\[
\chi_G(-x)=(-1)^{n}\chi_G(x)\,.
\]
\end{lemma}
\begin{proof}
One can mimic a standard proof for graphs (as in \cite[p.
11]{Big}; this result first appeared in a Chemistry paper
\cite{CoR}), or simply note that if one changes the signs of all
edges incident with vertices in one part then $\chi_G$ is
unchanged, yet every edge has then changed sign so that
$\chi_G(x)$ is changed to $(-1)^{n}\chi_G(-x)$.
\end{proof}

It will be convenient to extend the definition of bipartite to
cover any charged signed graph such that changing the sign of
every edge and charge produces a graph that is strongly
equivalent to the original. For (neutral) signed graphs, this
captures the usual definition of being bipartite. The extension of
Lemma \ref{L:bipartite} holds true for this larger class of
bipartite charged signed graphs, with the same proof.

A \emph{cycle} of length $r$ in a charged signed graph $G$ is a
list of distinct vertices $v_1, \ldots, v_r$ such that there is an
edge in $G$ between $v_i$ and $v_{i+1}$ ($1\le i<r$) and between
$v_1$ and $v_r$.  A charged signed graph without cycles is called
a (charged signed) \emph{forest}. A connected forest is called a
\emph{tree}.

\begin{lemma}[{\cite[Theorem 2.2]{CST}}]
Any charged signed forest is equivalent to one for which all the
edges are positive.
\end{lemma}
\begin{proof}
An easy induction on the number of vertices:
for the inductive step consider removing a leaf (a
vertex with exactly one neighbour), unless there are no edges.
\end{proof}

For detecting noncyclotomic integer symmetric matrices, the
following trivial and obvious sufficient condition can be useful.

\begin{lemma}\label{L:non-cyc}
Let $A$ be an $n\times n$ integer symmetric matrix.  If either
$\chi_A(2)<0$ or $(-1)^{n}\chi_A(-2)<0$, then $A$ is not
cyclotomic.
\end{lemma}


\begin{lemma}\label{L:smallcases}
Up to equivalence, the only indecomposable $1$-by-$1$ or
$2$-by-$2$ cyclotomic matrices  are
\[
(0)\,, (1)\,, (2)\,,
\left(\begin{array}{cc}0&1\\1&0\end{array}\right)\,,
\left(\begin{array}{cc}1&1\\1&0\end{array}\right)\,,
\left(\begin{array}{cc}1&1\\1&1\end{array}\right)\,,
\left(\begin{array}{cc}1&1\\1&-1\end{array}\right)\,
\text{and}\, \left(\begin{array}{cc}0&2\\2&0\end{array}\right)\,.
\]

Of these, the only maximal ones are $(2)$ and
$\left(\begin{array}{cc}0&2\\2&0\end{array}\right)$.
\end{lemma}
\begin{proof}
This is an easy computation, using Lemma \ref{L:non-cyc} to
constrain the matrix entries.
For example, to show that
$\left(\begin{array}{cc}0&2\\2&0\end{array}\right)$
is maximal,
suppose that
\[
A=\left(\begin{array}{ccc}0&2&a\\2&0&b\\a&b&c\end{array}\right)
\]
is cyclotomic.
To achieve $\chi(2)\ge0$ and $\chi(-2)\le 0$
requires both $-2(a+b)^2\ge0$ and $2(b-a)^2\le0$,
giving $a=b=0$, so that $A$ is not indecomposable.
\end{proof}

\begin{lemma}\label{L:reducetochargedsignedgraphs}
Apart from matrices equivalent to either $(2)$ or
$\left(\begin{array}{cc}0&2\\2&0\end{array}\right)$, any
indecomposable cyclotomic  matrix has all entries from the set
$\{0,1,-1\}$. In other words, it is the adjacency matrix of a
cyclotomic charged signed graph.
\end{lemma}
\begin{proof}
Let $A=(a_{ij})$ be an indecomposable cyclotomic matrix, not
equivalent to either $(2)$ or
$\left(\begin{array}{cc}0&2\\2&0\end{array}\right)$. Suppose first
that some diagonal entry of $A$ had modulus at least 2, say
$|a_{ii}|\ge2$. By interlacing (Lemma \ref{L:interlacing}), the
1-by-1 matrix $(a_{ii})$ is cyclotomic, and then by Lemma
\ref{L:smallcases} it equals $\pm (2)$ and is maximal, so equals
$A$, giving a contradiction.

Next suppose that some off-diagonal entry $a_{ij}$ had modulus at
least $2$. By interlacing, the 2-by-2 matrix
$\left(\begin{array}{cc}a_{ii}&a_{ij}\\a_{ij}&a_{jj}\end{array}\right)$
is cyclotomic, and by Lemma \ref{L:smallcases} this must equal
$\pm\left(\begin{array}{cc}0&2\\2&0\end{array}\right)$, and is
maximal, so equals $A$. Again we have a contradiction.

Thus no entry of $A$ has modulus greater than 1.
\end{proof}

We conclude that, apart from two (up to equivalence) trivial
examples, all indecomposable cyclotomic matrices are the adjacency
matrices of connected cyclotomic charged signed graphs. Thus
Theorem \ref{T-ISM} follows from Theorems \ref{T-CSG} and
\ref{T-LCSG}, and  we can restrict our attention to charged signed
graphs.

\section{Representation via Gram matrices}
\subsection{Gram matrices and line systems}
Let $A$ be the adjacency matrix of a cyclotomic charged signed graph
with $n$ vertices. In particular, $A$ has all eigenvalues at least
$-2$. Hence $A+2I$ is positive semi-definite. This implies that we
can find vectors ${\bf w}_1$, \ldots, ${\bf w}_n$ in real
$n$-dimensional space such that $A+2I$ is their Gram matrix: the
$(i,j)$-entry of $A+2I$ is the dot product of ${\bf w}_i$ and ${\bf
w}_j$. The dimension of the space spanned by the ${\bf w}_i$ might
of course be smaller than $n$.

A particularly simple case is that of a signed graph, where there
are no charges. Then the diagonal entries of $A+2I$ all equal 2, so
that the vectors ${\bf w}_i$ all have length $\sqrt{2}$. Moreover
the lines spanned by the ${\bf w}_i$ meet each other with angles
$\pi/3$ or $\pi/2$. In the language of \cite{CvL} we have
represented our signed graph in a line system, and if the graph is
connected then the line system is indecomposable.
If we change the
sign of one of our Gram vectors, then the line that it spans is
unchanged, and the new Gram matrix is equivalent to the old one: we
have just changed the sign of all edges incident with the vertex
that corresponds to our Gram vector. Since we are working up to
equivalence, we can fix (at our discretion) the direction of each
line in our system.

Indecomposable line systems have been classified. Every such line
system is contained in a maximal one. It follows that every
cyclotomic connected signed graph is contained in a maximal one.
Moreover we can hunt for these by looking inside the maximal
indecomposable line systems. These are $\mathcal D_n$ ($n\ge4$)
and $\mathcal E_8$,
which we now describe.

\subsection{The line system and signed graph $\mathcal D_n$}\label{SS:Dn}
Fix $n\ge 2$, and let ${\bf e}_1$, \ldots, ${\bf e}_n$ be an
orthonormal basis for ${\mathbb R}^n$. The signed graph $\mathcal
D_n$ has $n(n-1)$ vertices, represented by the vectors
\[
{\bf e}_i \pm {\bf e}_j\qquad (1\le i < j \le n)\,.
\]
Adjacency of unequal vertices is given by the dot product of the
corresponding vectors, which always equals one of 0, 1, $-1$.  If
$A$ is the adjacency matrix of $\mathcal D_n$, then $A+2I$ is the
Gram matrix of the set of vectors.

\subsection{The line system and signed graph $\mathcal E_8$}\label{SS:E8}
Let ${\bf e}_1$, \ldots, ${\bf e}_8$ be an orthogonal basis for
${\mathbb R}^8$, where, in contrast to the previous subsection,
each ${\bf e}_i$ has length $\sqrt{2}$. The signed graph $\mathcal
E_8$ has 120 vertices, represented by the vectors ${\bf e}_1$,
\ldots, ${\bf e}_8$ and 112 vectors of the form
\[
{\tfrac12}({\bf e}_i \pm {\bf e}_j \pm {\bf e}_k \pm {\bf e}_\ell
)\,,
\]
where $ijk\ell$ is one of the 14 strings
\[
\begin{array}{ccccccc}
1234\,,& 1256\,,& 1278\,,& 1357\,,& 1368\,,& 1458\,,& 1467\,, \\
2358\,,& 2367\,,& 2457\,,& 2468\,,& 3456\,,& 3478\,,& 5678\,.
\end{array}
\]
(The referee has pointed out that these strings are the supports
of the nontrivial words in the extended binary Hamming code of
length $8$.)

As for $\mathcal D_n$, adjacency of unequal vertices is given by
the dot product (one of 0, 1, $-1$).

As a notational convenience, the vertices of $\mathcal E_8$ will
be written as strings of digits, some of them overlined. Single
digits $1, \ldots, 8$ refer to the basis vectors ${\bf e}_1$,
\ldots, ${\bf e}_8$. Strings of four digits, with any of the last
three overlined, refer to the vectors $({\bf e}_i \pm {\bf e}_j
\pm {\bf e}_k \pm {\bf e}_\ell )/2$, with overlining indicating a
minus sign. For example, $14\bar{6}\bar{7}$ indicates the vector
$({\bf e}_1 + {\bf e}_4 - {\bf e}_6 - {\bf e}_7)/2$.

We sum up this discussion with the following result. For the proof
one trivially adapts to signed graphs the argument for graphs in
Chapter 3 of \cite{CvL}, noting that the fact that we can have
negative edges makes the argument significantly easier.

\begin{proposition}\label{T:DnorE8}
Up to equivalence, the only (neutral) connected signed graphs that
have all their eigenvalues in $[-2,\infty)$ are the connected
subgraphs of $\mathcal D_n$ ($n\ge 2$) and of $\mathcal E_8$.
\end{proposition}

Signed graphs with all their eigenvalues in $[-2,\infty)$ have
been studied earlier by Vijayakumar \cite{V}, Singhi  and
Vijayakumar \cite{VS} and Ray-Chaudhuri,  Singhi and Vijayakumar
\cite{RSV}.

\section{Cyclotomic signed graphs}
In this section we prove Theorem \ref{T-CSG}, and so classify all
cyclotomic signed graphs. The plan is as follows. First we find
all the connected cyclotomic signed graphs that contain triangles
(triples of vertices with each pair being adjacent). Then, in view
of Proposition \ref{T:DnorE8}, it suffices to consider
triangle-free subgraphs of $\mathcal D_n$ and $\mathcal E_8$. We
find all maximal triangle-free subgraphs of $\mathcal D_n$, and
observe the remarkable fact that they are all cyclotomic. We then
find all maximal triangle-free subgraphs of $\mathcal E_8$: these
are not all cyclotomic, and so we need to search among their
subgraphs for any new maximal connected cyclotomic signed graphs
that had not already been found as subgraphs of some $\mathcal
D_n$.

\subsection{Reduction to triangle-free graphs}
\begin{lemma}\label{L:addtotriangle}
Suppose that $G$ is a cyclotomic signed graph
that contains a triangle on vertices $v$, $w$, $x$
(the signs of the three edges being arbitrary).
If $z$ is a fourth vertex in $G$ then $z$ is
a neighbour of an even number of $v$, $w$, $x$.
\end{lemma}
\begin{proof}
Direct computation of the small number of cases. One finds that if
$z$ is a neighbour of one or three of $v$, $w$, $x$ then the
subgraph induced by $v$, $w$, $x$, $z$ is not cyclotomic,
contradicting $G$ being cyclotomic, by interlacing.
\end{proof}

If $z$ is a neighbour of exactly two of $v$, $w$, $x$, then the
subgraph induced by $v$, $w$, $x$, $z$ is not always cyclotomic, and
the next lemma describes the extra condition on the signs of the
edges that is required for a cyclotomic graph.
\begin{lemma}\label{L:parity}
If $G$ is a cyclotomic signed graph containing two triangles that
share an edge, then one triangle has an even number of negative
edges, and the other has an odd number of negative edges.
\end{lemma}
\begin{proof}
If two triangles share an edge and the parities of
the numbers of negative edges in the two triangles
are equal, then one quickly checks that a suitable equivalence
will make all the edges on both triangles positive.
But then the subgraph induced by the two triangles is
not cyclotomic (it has $(1+\sqrt{17})/2$
as an eigenvalue), and by interlacing neither is $G$.
\end{proof}

\begin{corollary}\label{C:triplysharededge}
If $G$ is a cyclotomic signed graph, then no
three triangles can share a single edge.
\end{corollary}

\begin{corollary}\label{C:tetrahedron}
If $G$ is a cyclotomic signed graph, then
it does not contain a tetrahedron as an induced
subgraph.
\end{corollary}
This latter Corollary also follows from Lemma
\ref{L:addtotriangle}.

\begin{lemma}\label{L:triangle}
If $G$ is a connected cyclotomic signed graph that contains a
triangle, then it is equivalent to a subgraph of the signed
octahedron $T_6$ of Figure \ref{F-oct}.
\end{lemma}

\begin{proof}
Suppose that $G$ is a connected cyclotomic signed graph that
contains a triangle, on vertices $v_1$, $v_2$, $v_3$. By a suitable
equivalence, we may suppose that the three edges of this triangle
are all positive. If $G$ contains no other vertices then we are
done.

Otherwise suppose that $v_4$ is another vertex of $G$, joined to
$v_1$, say. By Lemma \ref{L:addtotriangle}, $v_4$ is adjacent to
exactly one other of the $v_i$. Relabelling if necessary, we suppose
that $v_4$ is adjacent to $v_1$ and $v_2$. If $G$ contains no other
vertices then we are done.

Otherwise $G$ contains a fifth vertex $v_5$, adjacent to at least
one of $v_1$, $v_2$, $v_3$, $v_4$. By Lemma \ref{L:addtotriangle},
$v_5$ is adjacent to two vertices on one of the triangles
$v_1v_2v_3$, $v_1v_2v_4$, and hence is adjacent to one of $v_1$ or
$v_2$. By Corollary \ref{C:triplysharededge}, $v_5$ cannot be
adjacent to both $v_1$ and $v_2$. Without loss of generality, $v_5$
is adjacent to $v_1$. By Lemma \ref{L:addtotriangle} (using
triangles $v_1v_2v_3$ and $v_1v_2v_4$), $v_5$ is also adjacent to
both $v_3$ and $v_4$.
If $G$ contains no other vertices then we are done.

Otherwise $G$ contains a sixth vertex, $v_6$, adjacent to one of
$v_2$, $v_3$, $v_4$, $v_5$ (it cannot be adjacent to $v_1$, or else
by Lemma \ref{L:addtotriangle} it would be adjacent to one of the
others, producing three triangles sharing an edge, contrary to
Corollary \ref{C:triplysharededge}). Applying Lemma
\ref{L:addtotriangle} repeatedly, we see that $v_6$ must be adjacent
to all of $v_2$, $v_3$, $v_4$, $v_5$.

We now have a subgraph of $G$ that is equivalent to the signed
octahedron pictured in the Lemma (by Lemma \ref{L:parity} the parity
of the number of negative edges on faces sharing an edge must
differ, and up to equivalence one sees that there is just one choice
of signs).

Finally, $G$ can have no more vertices, as each existing triangle
shares each of its edges with another: we cannot adjoin a new vertex
in a way that is compatible with both Lemma \ref{L:addtotriangle}
and Corollary \ref{C:triplysharededge}.
\end{proof}

\begin{corollary}\label{C:degreebound}
In a cyclotomic signed graph $G$, each vertex has degree at most
$4$.

\end{corollary}
\begin{proof}
If $G$ contains a triangle then it is equivalent to
a subgraph of the signed octahedron, and hence has
maximal degree at most 4.
We may therefore assume that $G$ is triangle-free.

If $G$ has a vertex $v$ of degree at least 5, then $v$ has
neighbours $v_1$, \ldots, $v_5$ say (and possibly others), and since
$G$ is triangle-free there are no edges between any pair of $v_1$,
\ldots, $v_5$. By computation the starlike subgraph induced by $v$,
$v_1$, \ldots, $v_5$ is not cyclotomic (up to equivalence all the
edges are positive, so there is only one case to compute). This
contradicts $G$ being cyclotomic, by interlacing.
\end{proof}

\subsection{The maximal triangle-free subgraphs of $\mathcal D_n$}
After Lemma \ref{L:triangle}, our search for connected cyclotomic
signed graphs can be restricted to triangle-free connected
cyclotomic signed graphs. After Proposition \ref{T:DnorE8} we can
hunt for these triangle-frees as subgraphs of one of the $\mathcal
D_n$, or of $\mathcal E_8$. Here we deal with the $\mathcal D_n$,
classifying all the maximal triangle-free subgraphs. Fortunately
for us (in view of our ultimate goal) these subgraphs are all
cyclotomic.

For $v= {\bf e}_i  \pm {\bf e}_j\in \mathcal D_n$ --- so that
$i<j$
--- define the \textit{conjugate} vertex $v^*$ to be ${\bf e}_i
\mp {\bf e}_j$. If $v= {\bf e}_i  \pm {\bf e}_j$, then we say that
$v$ \textit{includes} ${\bf e}_i$ and ${\bf e}_j$. Note that $v$
and $v^*$ have the same neighbours in $\mathcal D_n$.

\begin{lemma}\label{L-conjugates}
Let $G$ be a maximal triangle-free subgraph of $\mathcal D_n$. If
$v$ is a vertex of $G$, then so is $v^*$.
\end{lemma}

\begin{proof}
If $v^* a b$ is a triangle in $G$, then so is $vab$. Hence if $G$
contained $v$ but not $v^*$ we could add $v^*$ to the vertex set and
get a larger triangle-free signed graph, contradicting the
maximality of $G$.
\end{proof}

\begin{lemma}\label{L-atmost4}
Let $G$ be a maximal triangle-free subgraph of $\mathcal D_n$.
Each ${\bf e}_i$ is included in at most four vertices of $G$.
\end{lemma}

\begin{proof}
If ${\bf e}_i$ is included at all, then take a vertex $v$
including  ${\bf e}_i$ and ${\bf e}_j$ ($j\ne i$).

Suppose first that there exists a vertex
 $w$ in $G$ that includes ${\bf e}_i$ and ${\bf e}_k$ for some
other $k$ ($k\ne i$, $k\ne j$). Then if $x$ is a vertex of $G$ that
includes ${\bf e}_i$ and ${\bf e}_\ell$
 ($\ell\ne i$)
we must have either $\ell=j$ or $\ell=k$, or else $vwx$ would be a
triangle. Hence ${\bf e}_i$ is included exactly four times, in $v$,
$w$, $v^*$, $w^*$.

If no such $w$
exists, then ${\bf e}_i$ is included in exactly two vertices,
$v$ and $v^*$.
\end{proof}

\begin{lemma}\label{L-exactly4}
Let $G$ be a maximal triangle-free subgraph of $\mathcal D_n$. The
maximum degree of $G$ is at most 4. Moreover if a vertex $v$ in
$G$ has distinct neighbours $a$ and $b$ with $a\ne b^*$, then $v$
has four neighbours, $a$, $b$, $a^*$, $b^*$.
\end{lemma}

\begin{proof}
Take  any vertex $v$ in $G$. By relabelling, and moving to $v^*$
if necessary, we can suppose that $v={\bf e}_1+{\bf e}_2$. Let $w$
be a neighbour of $v$. Again after relabelling, and so on, we can
suppose that $w={\bf e}_2+{\bf e}_3$. Then $w^*$ is also a
neighbour of $v$. If $v$ has a third neighbour $x$, then, in the
same way, we can suppose that $x={\bf e}_1+{\bf e}_4$.  Note that
$x$ cannot include ${\bf e}_2$, by Lemma \ref{L-atmost4}. Then
$x^*$ is a fourth neighbour of $v$. By  Lemma \ref{L-atmost4}
again, there can be no more neighbours, as these would have to
include either ${\bf e}_1$ or ${\bf e}_2$, both of which have been
included four times already (in $v$, $v^*$, $x$, $x^*$ and in $v$,
$v^*$, $w$, $w^*$ respectively).
\end{proof}

Recall that a path $v_1v_2\dots v_m$ in $G$ is a sequence of
distinct vertices $v_i$ in $G$ with $v_i$ adjacent to $v_{i+1}$
for $i=1,\dots,m-1$.

\begin{lemma}\label{L-maxpathiscycle}
Let $G$ be a maximal connected triangle-free subgraph of $\mathcal
D_n$, where $n\ge 4$. Let $P=v_1v_2\ldots v_m$ be a path in $G$,
maximal subject to no $v_i$ equalling any $v_j^*$. Then
\begin{itemize}
 \item $v_1$ and $v_m$ are adjacent, so that the induced subgraph on
 the vertices of $P$ is a cycle.
 \item
$P^*:=v_1^* \ldots v_m^*$ is a path in $G$ disjoint from $P$, and
$G$ is the subgraph spanned by $P$ and $P^*$.
\end{itemize}
\end{lemma}

\begin{proof}
First suppose that $v_1$ and $v_m$ are not adjacent. By Lemma
\ref{L-conjugates}, $P^*$ is a subgraph of $G$.  No vertex in $P$
can have more than two neighbours in $P$, else together with its
neighbours in $P^*$ it would have more than four neighbours in
$G$, contradicting Lemma \ref{L-exactly4}. Without loss of
generality, $v_1 = {\bf e}_1 + {\bf e}_2$, $v_2 = {\bf e}_2 + {\bf
e}_3$, \ldots, $v_{m-1} = {\bf e}_{m-1}+{\bf e}_m$, $v_m = {\bf
e}_m + {\bf e}_{m+1}$.

Now for $2\le i \le m-1$, $v_i$ has neighbours $v_{i-1}$,
$v_{i+1}$, $v_{i-1}^*$, $v_{i+1}^*$, so has no other neighbours in
$G$, by Lemma \ref{L-exactly4}.  By maximality of $P$, $v_1$ and
$v_m$ have no neighbours in $G$ that are not in $P$ or $P^{*}$, so
$P$ and $P^*$ span a component of $G$, and hence span $G$.  But
then we could add ${\bf e}_1 + {\bf e}_{m+1}$ to $G$ without
introducing triangles, contradicting maximality of $G$.

Thus $v_1$ and $v_m$ are adjacent, and without loss of generality
$v_1 = {\bf e}_1 + {\bf e}_2$, $v_2 = {\bf e}_2 + {\bf e}_3$,
\ldots, $v_{m-1} = {\bf e}_{m-1}+{\bf e}_m$, $v_m = {\bf e}_1 +
{\bf e}_{m}$.  Now each element of $P\cup P^*$ has four neighbours
in $P\cup P^*$, so no others, and again $P$ and $P^*$ span the
whole of $G$.
\end{proof}

The proof of Lemma \ref{L-maxpathiscycle} establishes the first
sentence of the next result.

\begin{proposition}\label{justabove}
Every maximal connected triangle-free signed graph that is a
subgraph of some $\mathcal D_n$ ($n \ge 4$) is equivalent to one
with vertex set of the form
\[
{\bf e}_1 + {\bf e}_2, {\bf e}_2+{\bf e}_3, \ldots, {\bf
e}_{m-1}+{\bf e}_m,
 {\bf e}_1+{\bf e}_m,
{\bf e}_1 - {\bf e}_2, {\bf e}_2 - {\bf e}_3, \ldots, {\bf
e}_{m-1}-{\bf e}_m,
 {\bf e}_1-{\bf e}_m\,,
\]
 for some $m$ in the range $4\le m \le n$.   Moreover, every
such graph is cyclotomic, and is a maximal connected cyclotomic
signed graph.
\end{proposition}

\begin{proof}
It remains to prove that such a graph $G$ is cyclotomic
(maximality as a connected cyclotomic signed graph follows from
Corollary \ref{C:degreebound}). For $n=4$ one gets this by
computation (or an easy adaptation of the following argument). For
$n>4$, note that if $v$ and $w\ne v^*$ are distance 2 apart in
$G$, then there are exactly two 2-paths from $v$ to $w$, one along
edges of the same sign, and one along edges of opposite sign.
There are four 2-paths from $v$ to $v^*$, two along edges of the
same sign, and two along edges of opposite sign. Hence (with $A$
the adjacency matrix of $G$) all off-diagonal entries of $A^2$ are
zero. Since each vertex has degree 4, we deduce that $A^2 = 4I$.
Hence all the eigenvalues are either $2$ or $-2$, so $G$ is
cyclotomic.
\end{proof}

A nice representative of the equivalence class of the maximal
cyclotomic signed graph, denoted $T_{2n}$ in the Theorem, described in Proposition \ref{justabove} is
obtained by replacing the vertex ${\bf e}_1-{\bf e}_n$ by ${\bf
e}_n - {\bf e}_1$. Then one of the $n$-cycles (say $v_1 v_2 \cdots
v_n$) has all positive edges, and the other ($v_1^* v_2^* \cdots
v_n^*$) has all negative edges. The linking edges of the form $v_i
v_{i+1}^*$ (interpreted cyclically) are all positive, and those of
the form $v_i v_{i-1}^*$ are all negative. One gets a nice picture
if the two cycles are viewed as the ends of a cylinder.
Alternatively, the graph can be drawn on a torus without
crossings, wrapping each cycle round the torus in such a way that
it cannot be shrunk to a point (as in Figure \ref{F-torus}).

\subsection{The maximal triangle-free subgraphs of $\mathcal E_8$}\label{S:E8tfree}
The search for triangle-free subgraphs of $\mathcal E_8$ (up to
equivalence) was done by computer, using moderately intelligent
backtracking. A lexicographical ordering was given to the 120
vertices, and a set of equivalences of $\mathcal E_8$ was
precomputed (each as an explicit permutation of the 120 vertices),
as follows.

We can change the sign of any ${\bf e}_i$: for any $i\in\{1, 2, 3,
4, 5, 6, 7, 8\}$, we can swap the roles of $i$ and $\bar{i}$,
which preserves all dot products. Then flip the sign of any vector
that is no longer a vertex of $\mathcal E_8$ to induce an
equivalence of $\mathcal E_8$. If $G$ is a signed subgraph of
$\mathcal E_8$, then applying this process gives an equivalent
(but perhaps different) subgraph.

Some, but not all, permutations of $\{1, 2, 3, 4, 5, 6, 7, 8\}$
induce a permutation of the lines spanned by the vertices of
$\mathcal E_8$ (and hence induce a permutation of the vertices of
$\mathcal E_8$). For any string $ijk\ell$ that appears as a
vertex, we can apply elements of the Klein 4-group acting on
$\{i,j,k,\ell\}$ to induce a permutation of the vertices of
$\mathcal E_8$. For example, if we apply $(12)(56)$ to the vertex
$1\bar{2}3\bar{4}$ we get the vector $\bar{1}23\bar{4}$, which
spans the same line as the vertex $1\bar{2}\bar{3}4$, so the image
of $1\bar{2}3\bar{4}$ under $(12)(56)$ is $1\bar{2}\bar{3}4$. Note
that such a transformation might not be an isomorphism of signed
graphs (since some of the vertices may be switched)
but will be an equivalence. Again, applying this to a subgraph of
$\mathcal E_8$ will give an equivalent subgraph.

Also, if $ijk\ell$ is a vertex of $\mathcal E_8$, then we can
perform a change of basis by the following four swaps:
$i\leftrightarrow ijk\ell$, $j\leftrightarrow
ij\bar{k}\bar{\ell}$, $k\leftrightarrow i\bar{j}k\bar{\ell}$,
$l\leftrightarrow i\bar{j}\bar{k}\ell$. This induces an
equivalence on $E$. (It is enough to check that this works for
$ijk\ell=1234$, and then use the previous symmetries to reduce to
this case.)

Starting with $S$ being empty, the search grew $S$ by adding the
smallest possible vertices (with respect to the chosen ordering)
whilst (i) maintaining triangle-freeness, and (ii) checking that
none of the above equivalences of $\mathcal E_8$ would map the
enlarged $S$ to a lexicographically earlier set. The use of
equivalences was hugely powerful in cutting down on the number of
sets $S$ considered by rejecting most sets at an early stage. When
no more vertices could be added, the set $S$ was tested for
maximality, and maximal triangle-frees were written to a file.
Then backtracking was done to find the next candidate for $S$.

The following twenty inequivalent maximal triangle-free subgraphs
were found.
\begin{itemize}

\item[{\bf G1}] $1, 2,  3, 4, 5, 6, 7, 8, 1234,
12\bar{3}\bar{4}, 1\bar{2}3\bar{4}, 1\bar{2}\bar{3}4,
5678, 56\bar{7}\bar{8}, 5\bar{6}7\bar{8}, 5\bar{6}\bar{7}8.$

This is two copies of the toral tesselation $T_8$.

\item[{\bf G2}] $1, 2,  3, 4, 5, 6, 7, 8, 1234,
12\bar{3}\bar{4}, 1\bar{2}5\bar{6}, 1\bar{2}\bar{5}6,
3\bar{4}56, 3\bar{4}\bar{5}\bar{6}.$

This comprises two isolated vertices plus
$T_{12}$.

\item[{\bf G3}] $1, 2,  3, 4, 5, 6, 7, 8, 1234,
12\bar{3}\bar{4}, 1\bar{2}5\bar{6}, 1\bar{2}\bar{5}6,
3\bar{4}7\bar{8}, 3\bar{4}\bar{7}8, 5678, 56\bar{7}\bar{8}.$

This is $T_{16}$.

\item[{\bf G4}] $1, 2,  3, 4, 5, 6, 7, 8, 1234, 1\bar{2}5\bar{6},
1\bar{3}\bar{5}7, 1\bar{4}6\bar{7}, 2\bar{3}5\bar{8},
2\bar{4}\bar{6}8, 3\bar{4}7\bar{8}, 5678.$

This is the hypercube $S_{16}$. 

\item[{\bf G5}] $1, 2,  3, 4, 5, 6, 7, 8, 1234,
1\bar{2}5\bar{6}, 1\bar{3}\bar{5}7, 1\bar{4}6\bar{7},
2\bar{3}\bar{6}\bar{7}, 2\bar{4}57, 3\bar{4}\bar{5}\bar{6}.$

This is an isolated vertex plus $S_{14}$.

\item[{\bf G6}] $1, 2,  3, 4, 5, 6, 1234,
12\bar{3}\bar{4}, 1\bar{2}5\bar{6}, 1\bar{2}\bar{5}6,
3\bar{4}7\bar{8}, 3\bar{4}\bar{7}8, 567\bar{8}, 56\bar{7}8.$

This is $T_{14}$.

\item[{\bf G7}] $1, 2,  3, 4, 5, 6, 1234,
12\bar{3}\bar{4}, 1\bar{2}7\bar{8}, 1\bar{2}\bar{7}8,
3\bar{4}7\bar{8}, 3\bar{4}\bar{7}8, 5678, 56\bar{7}\bar{8}.$

This is a square plus  $T_{10}$.

\item[{\bf G8}] $1, 2,  3, 5, 1278, 14\bar{6}\bar{7},
2\bar{4}6\bar{8}, 3456, 3\bar{4}\bar{7}8,
 5\bar{6}7\bar{8}.$

10 vertices, 2 cyclotomic components (both are 5-cycles).
\end{itemize}

By Corollary \ref{C:degreebound}, $S_{14}$ and $S_{16}$ are
maximal.

In the remaining cases, the larger component was noncyclotomic.

\begin{itemize}

\item[{\bf G9}] $1, 2,  3, 4, 5, 6, 1234, 1\bar{2}5\bar{6},
1\bar{3}\bar{5}7, 1\bar{4}6\bar{7}, 3\bar{4}7\bar{8}, 567\bar{8}.$

12 vertices, 1 component, 29 maximal cyclotomic subgraphs
(maximal in the sense that no larger subgraph of {\bf G9} is
cyclotomic).

\item[{\bf G10}] $1, 2,  3, 4, 5, 6, 1234, 1\bar{2}7\bar{8},
1\bar{3}5\bar{7}, 1\bar{4}\bar{5}8, 3\bar{4}\bar{7}\bar{8}, 5678.$

12 vertices, 1 component, 13 maximal cyclotomic subgraphs.

\item[{\bf G11}] $1, 2,  3, 4, 5, 6, 1234, 1\bar{2}7\bar{8},
1\bar{3}5\bar{7}, 2\bar{4}57, 3\bar{4}\bar{7}8.$

11 vertices, 2 components (one being a single vertex), 15 maximal cyclotomic subgraphs.

\item[{\bf G12}] $1, 2,  3, 4, 5, 6, 1234, 1\bar{2}7\bar{8},
1\bar{3}5\bar{7}, 2\bar{4}6\bar{8}, 3\bar{4}\bar{7}8, 5678.$

12 vertices, 1 component, 15 maximal cyclotomic subgraphs.

\item[{\bf G13}] $1, 2,  3, 4, 5, 6, 1234, 1\bar{2}7\bar{8},
1\bar{3}5\bar{7}, 2\bar{4}6\bar{8}, 5678, 5\bar{6}7\bar{8}.$

12 vertices, 1 component, 19 maximal cyclotomic subgraphs.

\item[{\bf G14}] $1, 2,  3, 4, 5, 6, 1278, 1\bar{2}7\bar{8},
13\bar{5}\bar{7}, 235\bar{8}, 3478, 5\bar{6}\bar{7}8.$

12 vertices, 1 component, 17 maximal cyclotomic subgraphs.

\item[{\bf G15}] $1, 2,  3, 4, 5, 6, 1278, 1\bar{2}7\bar{8},
13\bar{5}\bar{7}, 24\bar{6}\bar{8}, 3478, 5\bar{6}\bar{7}8.$

12 vertices, 1 component, 37 maximal cyclotomic subgraphs.

\item[{\bf G16}] $1, 2,  3, 4, 5, 1234, 1\bar{2}5\bar{6},
1\bar{3}6\bar{8}, 2\bar{4}\bar{6}8, 3\bar{4}7\bar{8}, 56\bar{7}8.$

11 vertices, 1 component, 44 maximal cyclotomic subgraphs.

\item[{\bf G17}] $1, 2,  3, 4, 5, 1234, 1\bar{2}7\bar{8},
1\bar{3}\bar{6}8, 2\bar{4}6\bar{8}, 3\bar{4}\bar{7}8, 5678.$

11 vertices, 2 components: $K_2$ plus a 9-vertex component; 36
maximal cyclotomic subgraphs.

\item[{\bf G18}] $1, 2,  3, 4, 5, 1256, 1\bar{2}7\bar{8},
13\bar{5}\bar{7}, 24\bar{6}\bar{8}, 3478, 3\bar{4}7\bar{8},
5\bar{6}\bar{7}\bar{8}.$

12 vertices, 1 component, 3-regular, 45 maximal
cyclotomic subgraphs.

\item[{\bf G19}] $1, 2,  3, 4, 5, 1278, 13\bar{6}\bar{8},
1\bar{4}6\bar{7}, 236\bar{7}, 2\bar{4}\bar{6}\bar{8}, 3\bar{4}78,
5\bar{6}\bar{7}8.$

12 vertices, 2 components: $K_2$ plus a 3-regular 10-vertex
component, equivalent to the Petersen graph
(switch at vertex $4$ to get it)
57 maximal cyclotomic
subgraphs.

\item[{\bf G20}] $1, 2,  3, 5, 1234, 1\bar{2}5\bar{6},
1\bar{3}6\bar{8}, 2\bar{4}5\bar{7}, 3\bar{4}\bar{7}\bar{8},
567\bar{8}.$

10 vertices, 1 component, 23 maximal
cyclotomic subgraphs.

\end{itemize}

For each of the noncyclotomic components listed above, it was
checked by computer that none of their cyclotomic subgraphs are
maximal cyclotomic graphs.

This completes the proof of Theorem \ref{T-CSG}.


\subsection{An alternative view of the cyclotomic subgraphs of $\mathcal E_8$}
Let $G$ be a cyclotomic subgraph of $\mathcal E_8$ that is not
equivalent to a subgraph of any $\mathcal D_r$, with vertices
given by vectors ${\bf v}_1$, \ldots, ${\bf v}_n$, contained in
8-dimensional real space. Then $-G$, obtained from $G$ by changing
the signs of all edges and charges, is also cyclotomic. Now $-G$
is equivalent to $G$, so cannot be represented in any line system
$\mathcal D_r$, so must be represented in the line system
$\mathcal E_8$, and hence the vertices of $-G$ can be represented
as vectors ${\bf w}_1$, \ldots, ${\bf w}_n$, where for each $i$
either ${\bf w}_i$ or $-{\bf w}_i$ is in the signed graph
$\mathcal E_8$.

We can view the concatenated vectors $[{\bf v}_1,{\bf w}_1]$,
\ldots, $[{\bf v}_n,{\bf w}_n]$ as elements of 16-dimensional real
space, a subset of the 28800 vectors $[{\bf v},{\bf w}]$ where
${\bf v}\in \mathcal E_8$, $\pm{\bf w}\in \mathcal E_8$. Moreover,
since the ${\bf w}_i$ represent $-G$, we have
\[
{\bf w}_i \cdot {\bf w}_j = - {\bf v}_i \cdot {\bf v}_j
\]
for all $i\ne j$.  This implies that
\[
[{\bf v}_i,{\bf w}_i]\cdot[{\bf v}_j,{\bf w}_j] = 0
\]
for all $i\ne j$: our concatenated vectors $[{\bf v}_1,{\bf w}_1]$,
\ldots, $[{\bf v}_n,{\bf w}_n]$ are pairwise orthogonal. Since these
vectors lie in 16-dimensional space, we must have $n\le 16$, as is
confirmed by the examples computed in Section \ref{S:E8tfree}.

Conversely, suppose that we take any orthogonal subset $[{\bf
v}_1,{\bf w}_1]$, \ldots, $[{\bf v}_n,{\bf w}_n]$ of the 28800
vectors considered above, with the constraint that ${\bf v}_1$,
\ldots, ${\bf v}_n$ are distinct. Then the signed graph $G$ with
vertices ${\bf v}_1$, \ldots, ${\bf v}_n$ (and adjacency of
unequal vertices given by the dot product) is cyclotomic, for both
$G$ and $-G$ are represented in the line system $\mathcal E_8$,
with ${\bf w}_1$, \ldots, ${\bf w}_n$ spanning the lines that
represent $-G$.

\subsection{Remark on maximal cyclotomic (unsigned) graphs}
\begin{figure}[h]
\begin{center}
\leavevmode
\psfragscanon
\psfrag{A}[l]{$\tilde A_n$}
\psfrag{D}[l]{$\tilde D_n$}
\psfrag{6}[l]{$\tilde E_6$}
\psfrag{7}{$\tilde E_7$}
\psfrag{8}{$\tilde E_8$}
\includegraphics[scale=0.4]{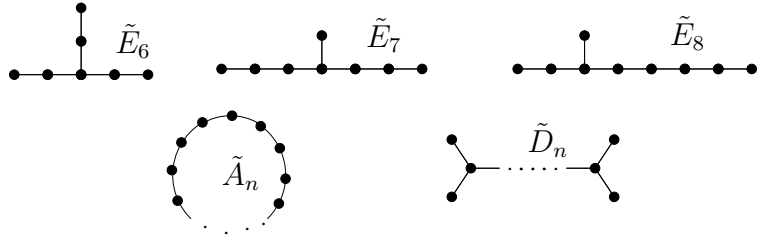}
\end{center}
\caption{The maximal connected cyclotomic graphs $\tilde E_6, \tilde E_7,
 \tilde E_8, \tilde A_n (n\geq 2)$ and $ \tilde
D_n (n\geq 4)$. The number of vertices is $1$ more than the index.
(From \cite{MS}).}\label{F-except}
\end{figure}

The maximal cyclotomic graphs classified by Smith are shown in
Figure \ref{F-except}. The $n$-cycle
$\tilde{A}_{n-1}$ and the graph $\tilde{D}_n$
 are subgraphs of $T_{2n}$, while the sporadic examples are all subgraphs
 of the hypercube $S_{16}$.  Unlike in the signed
case, however,  the maximal unsigned graphs are not visibly cyclotomic.

We can deduce Smith's classification as a corollary of Theorem
\ref{T-CSG}, by checking that these graphs are the only maximal
 (unsigned) subgraphs of the signed graphs of the Theorem.
A useful fact to use in this check is that the graphs $\tilde D_4$
and $\tilde D_5$, since they have $2$ as an eigenvalue, cannot be
a proper subgraph of any such graph. This is because otherwise the
graph would have an eigenvalue greater than $2$ --- see \cite[p.
4]{CvR}.

 We also note in passing that the classification of all graphs having all their
 eigenvalues in the open interval $(-2,2)$ follows from Smith's result. Such a
 graph is either a subgraph of $E_8$ or of some $D_n$ for $n\ge 8$ (Figure \ref{F-P}).
 Here, $E_8$ is $\tilde E_8$ (Figure \ref{F-except}) with its rightmost vertex removed (same as $U_5$ in
 Figure \ref{F-U}), and $D_n$ is $\tilde{D}_n$
 with a leaf removed. See also Theorem
\ref{T:nonnegative_matrix_open} below for a generalisation of this
result.

\section{Cyclotomic charged signed graphs}
We now embark upon the trickier task of proving Theorem
\ref{T-LCSG}, and so classifying all cyclotomic charged signed
graphs. The addition of charges means that we can no longer appeal
to Proposition \ref{T:DnorE8}, although the Gram matrix approach
will still prove extremely powerful.

\subsection{Excluded subgraphs I}\label{S:excludedI}
\begin{figure}[h]
\begin{center}
\leavevmode
\psfragscanon
\psfrag{X_1}{$X_1$}
\psfrag{X_2}{$X_2$}
\psfrag{X_3}{$X_3$}
\psfrag{X_4}{$X_4$}
\psfrag{X_5}{$X_5$}
\psfrag{X_6}{$X_6$}
\psfrag{X_7}{$X_7$}
\psfrag{X_8}{$X_8$}
\hbox{ \epsfxsize=3.3in
\epsffile{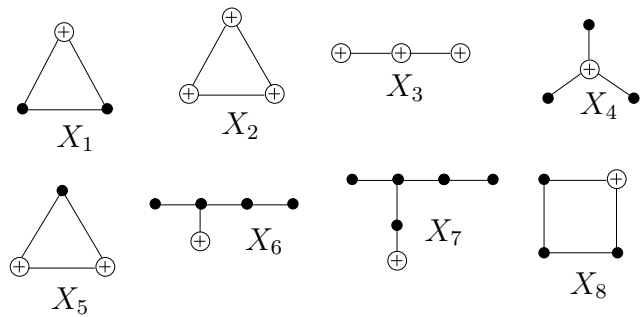} }
\end{center}
\caption{Excluded subgraphs I: some noncyclotomic charged signed graphs.}  \label{F-X}
\end{figure}

By interlacing, every subgraph of a cyclotomic charged signed
graph is cyclotomic. We can therefore exclude as subgraphs any
that are not cyclotomic. In particular, the following eight
non-cyclotomic charged signed graphs $X_1,\dots,X_8$ of Figure
\ref{F-X} (or anything equivalent to any of them) cannot be
subgraphs of any cylotomic charged signed graph.

\subsection{Excluded subgraphs II}\label{S:excludedII}
Certain cyclotomic charged signed graphs have the property that if
one tries to grow them to give larger connected cyclotomic graphs
then one always stays inside one of the maximal examples on the
following list: $S_7$, $S_8$, $S'_8$, $C^{++}_4$, $C^{+-}_4$, $C^{++}_6$,
$C^{+-}_6$, $T_6$. The process of proving that a cyclotomic graph
has this property is in principle simple, although perhaps
tedious, to carry out. Starting from the given graph, one
considers all possible ways of adding a vertex (up to equivalence)
such that the graph remains connected and cyclotomic. Check that
the resulting graphs are (equivalent to) subgraphs of one of
graphs on this list. Repeat with all the larger graphs found. If
the checks in this process are always valid, then, since the
process terminates, the original graph is suitable for exclusion.

\begin{figure}[h]
\begin{center}
\leavevmode
\psfragscanon
\psfrag{Y_1}{$Y_1$}
\psfrag{Y_2}{$Y_2$}
\psfrag{Y_3}{$Y_3$}
\psfrag{Y_4}{$Y_4$}
\psfrag{Y_5}{$Y_5$}\psfrag{Y_6}{$Y_6$}
\hbox{ \epsfxsize=4.3in
\epsffile{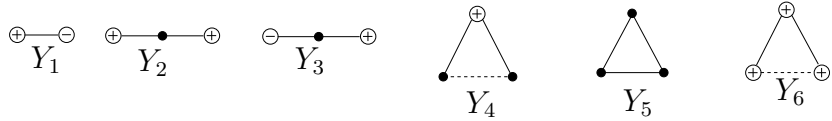} }
\end{center}
\caption{Excluded subgraphs II: Some cyclotomic charged signed
graphs that are contained as subgraphs of a maximal connected
cyclotomic charged signed graph only in one of the maximal graphs
$S_7$, $S_8$, $S'_8$, $C^{++}_4$, $C^{+-}_4$, $C^{++}_6$,
$C^{+-}_6$, $T_6$.} \label{F-Y}
\end{figure}

By this technique, the  six cyclotomic graphs $Y_1,\dots,Y_6$ of
Figure \ref{F-Y} (and anything equivalent to them) can be excluded
from future consideration.


\subsection{Charged and neutral components}
Let $G$ be a charged signed graph. We define the \emph{charged
subgraph} of $G$ to be the subgraph induced by all its charged
vertices, and the \emph{neutral subgraph} of $G$ to be the
subgraph induced by all its neutral vertices. The components of
the charged (respectively neutral) subgraph of $G$ will be called
the \emph{charged components} of $G$ (respectively the
\emph{neutral components} of $G$).

Our next task will be to show that the charged components of a
cyclotomic charged signed graph are tiny, provided that $G$ does
not contain $Y_1$, $Y_6$, or any equivalent subgraph.

\begin{lemma}\label{L:chargedcomponents}
Let $G$ be a cyclotomic charged signed graph that does not contain
any subgraph equivalent to $Y_1$ or $Y_6$ of Section
\ref{S:excludedII}. Then each charged component of $G$ contains at
most two vertices, necessarily of the same charge.
\end{lemma}
\begin{proof}
The last phrase is clear, since $Y_1$ is excluded as a subgraph.
Moreover the exclusion of $Y_1$ forces every charged component to
have all charges of the same sign, which by equivalence we may
assume to be all positive. Since graphs $X_2$ and $X_3$ of Section
\ref{S:excludedI} are not cyclotomic, and $Y_6$ is excluded by
assumption,  no charged component of $G$ can have as many as three
vertices.
\end{proof}

\subsection{Local geometric constraints}
\begin{lemma}\label{L:localgeometry1}
Let $G$ be a cyclotomic charged signed graph. Suppose that $G$
contains two nonadjacent neutral vertices $v$ and $w$ that have a
charged vertex $x$ as a common neighbour. Then $v$ and $w$ have
the same neighbours.
\end{lemma}
\begin{proof}
Adjacency being unchanged by equivalence, we may suppose that the
charge on $x$ has negative sign, and that the edges joining $v$
and $w$ to $x$ are positive. The subgraph induced by $v$, $w$, $x$
is then 
\noindent\phantom{ $x$ \text{is then}}  \psfragscanon
\psfrag{w}{$w$}
\psfrag{x}{$x$}
\psfrag{v}{$v$}\hbox{ \epsfxsize=0.8in
\epsffile{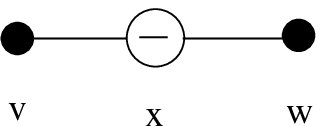} }
.


Since $G$ is cyclotomic, all eigenvalues of its adjacency matrix $A$
are in $[-2,\infty)$, so $A+2I$ is the Gram matrix of some set of
vectors. Let $\bf v$, $\bf w$, $\bf x$ be the Gram vectors
corresponding to $v$, $w$, $x$. Since $v$ and $w$ are neutral,
$\bf v$ and $\bf w$ have length $\sqrt{2}$. Since $x$ has a negative
charge, $\bf x$ has length 1. The angle between $\bf v$ and $\bf x$
is $\pi/4$, the angle between $\bf w$ and $\bf x$ is $\pi/4$, and
the angle between $\bf v$ and $\bf w$ is $\pi/2$. Hence $\bf v$,
$\bf w$, $\bf x$ are coplanar, with $\bf x$ in the direction of
${\bf v}+{\bf w}$. By consideration of their lengths we have
\begin{equation}\label{xvw}
2{\bf x}={\bf v} + {\bf w}\,.
\end{equation}

Now let $y$ be any other vertex of $G$, with
corresponding Gram vector $\bf y$.
Taking dot products with (\ref{xvw}) gives
\begin{equation}\label{ydotx}
2{\bf y}\cdot{\bf x} = {\bf y}\cdot{\bf v} + {\bf y}\cdot{\bf w}\,.
\end{equation}
The left hand side of (\ref{ydotx}) is an
even integer, hence the parities of the two integers on the
right must agree.
Hence $y$ is adjacent to $v$ if and only if it is adjacent
to $w$.
\end{proof}

\begin{lemma}\label{L:localgeometry2}
Let $G$ be a cyclotomic charged signed graph containing adjacent
charged vertices $v$ and $w$, where the signs on the charges for
$v$ and $w$ agree. Then $v$ and $w$ have the same neighbours.
\end{lemma}

\begin{proof}
Adjacency is preserved by equivalence, so we may suppose that the
charges on $v$ and $w$ are both negative, and that the edge
between $v$ and $w$ is positive. In the usual way, let $\bf v$
and $\bf w$ be Gram vectors corresponding to $v$, $w$. These have
length 1, and the angle
between them is zero, so $\bf v=\bf w$, although $v\ne w$. 
Hence $v$ and $w$ have the same neighbours.
\end{proof}

\subsection{Removing charged components I}
We now show that if a connected cyclotomic charged signed graph
does not contain a subgraph equivalent to any of the excluded
subgraphs of Section \ref{S:excludedII}, then it has a single
neutral component. As a first step, we show that certain charged
vertices can be deleted without disconnecting the graph.

\begin{lemma}\label{L:removecharge1}
Suppose that a connected cyclotomic charged signed graph $G$ has two
adjacent charged vertices $v$ and $w$, with the charges on $v$ and
$w$ having the same sign. Then the vertex $w$ can be deleted without
disconnecting $G$.
\end{lemma}
\begin{proof}
By Lemma \ref{L:localgeometry2} every neighbour of $v$ is a
neighbour of $w$ (and vice versa). Let $x$ and $y$ be any distinct
vertices in $G$, with neither of them being $w$. We must show that
there is a walk in $G$ from $x$ to $y$ that does not pass through
$w$. Certainly there is a path $v_1 v_2 \ldots v_r$ in $G$ from $x$
to $y$ ($v_1=x$, $v_r=y$). Suppose that this path contains $w$, say
$v_i = w$. If either $v_{i-1}$ or $v_{i+1}$ is $v$, then we can
simply remove $w$ from the path, since $v$ shares its neighbours.
Otherwise we can replace $w$ by $v$ in the path (producing a walk,
but perhaps no longer a path), again since $v$ and $w$ share their
neighbours.
\end{proof}

\begin{lemma}\label{L:removecharge2}
Let $G$ be a connected cyclotomic charged signed graph, with more
than three vertices, that contains no subgraph equivalent to either
$Y_1$ or $Y_4$ of Section \ref{S:excludedII}. Suppose that $G$ contains
two nonadjacent neutral vertices $v$ and $w$ that share a common
charged neighbour $x$ (as in Lemma \ref{L:localgeometry1}). Then $x$
can be deleted from $G$ without disconnecting the graph.
\end{lemma}
\begin{proof}
By the hypothesis on the number
of vertices in $G$, there is some
fourth vertex $y$ in $G$ that is adjacent
to one of $v$, $w$, $x$.

First we dispose of the cases where $y$ is adjacent to $x$. If $y$
has a charge, then since $Y_1$ is an excluded subgraph, $y$ and
$x$ have charges of the same sign. Then Lemma
\ref{L:removecharge1} shows that $x$ can be removed without
disconnecting $G$. If $y$ is neutral, then since subgraphs
equivalent to $Y_4$ are excluded, and a subgraph equivalent
to $X_1$ of Section \ref{S:excludedI} is impossible, $y$ cannot be
adjacent to either $v$ or $w$. But then $G$ would contain a
subgraph equivalent to $X_4$ of Section \ref{S:excludedI}, which
is not possible.

We may now suppose that $y$ is not adjacent to $x$, and more
strongly may suppose that $v$ and $w$ are the only neighbours of
$x$. By Lemma \ref{L:localgeometry1}, $v$ and $w$ share all their
neighbours. In particular, $y$ is adjacent to both $v$ and $w$.

Let $z_1$ and $z_2$ be any vertices in $G$
other than $x$.
It is enough to show that there is a
walk in $G$ from $z_1$ to $z_2$ that does
not pass through $x$.
Certainly there is a path $v_1 v_2 \ldots v_r$
from $z_1$ to $z_2$ ($v_1=z_1$, $v_r=z_2$).
Suppose that $x$ is on this path:
say $x = v_i$.
We know that $v_{i-1}$ and $v_{i+1}$
each equal one of $v$ and $w$.
We can therefore replace $x$ by $y$ in
our path to produce the desired walk.
\end{proof}

The requirement that $G$ has more than three
vertices is clearly necessary: if $v$, $w$, $x$
are the only vertices in $G$ then deleting
$x$ disconnects $G$.

\subsection{Removing charged components II}
\begin{lemma}\label{L:singleneutralcomponent}
Let $G$ be a connected cyclotomic charged signed graph that does
not contain a subgraph equivalent to $Y_1$, $Y_4$ or $Y_6$ of
Section \ref{S:excludedII}. Suppose further that $G$ has at least
four vertices. Then $G$ contains a single neutral component: all
charged vertices can be deleted without disconnecting $G$.
\end{lemma}

\begin{proof}
By Lemma \ref{L:chargedcomponents}, all charged components have at
most two vertices, and do not equal $Y_1$. By Lemma
\ref{L:removecharge1}, we can remove a charged vertex from any
charged component that has two vertices, without disconnecting $G$.
We are thus reduced to charged components containing only one
vertex.

If a charged vertex is a leaf, it can be removed without
disconnecting $G$.

If a charged vertex has two neutral neighbours, then since
subgraphs equivalent to $X_1$ and $Y_4$ are excluded we can
appeal to Lemma \ref{L:removecharge2} to see that this vertex can
be removed without disconnecting $G$.

No charged vertex can have three or more neutral neighbours, or
$G$ would contain a subgraph equivalent to one of $X_1$, $X_4$ or $Y_4$.
\end{proof}


\subsection{Growing the neutral component}
Let $G$ be a connected cyclotomic charged signed graph that
contains at least four vertices, at least one of which is charged,
but does not contain any of the excluded subgraphs $Y_1$, \ldots,
$Y_6$.  Then Lemma \ref{L:singleneutralcomponent} tells us that
$G$ has a single neutral component, $H$ say. By interlacing, $H$
is cyclotomic, and from the classification of all cyclotomic
signed graphs we know that $H$ is (equivalent to) a subgraph of
one of $\mathcal D_r$ (for some $r$), $S_{14}$ or $S_{16}$. We
treat each of these cases in turn.

\subsubsection{$H$ is equivalent to a subgraph of $\mathcal D_r$}
We may suppose that $r$ is minimal such that $\mathcal D_r$
contains a subgraph equivalent to $H$. Cases with $r\le 4$ can be
dealt with exhaustively by growing each possible $H$ in all
possible ways, adding only charged vertices, and checking that
each maximal connected cyclotomic  charged signed graph (maximal
subject to the neutral component being $H$) is contained in some
$C_{2k}^{++}$ or $C_{2k}^{+-}$. We may therefore suppose that
$r\ge 5$.

Working up to equivalence, we identify $H$ with some subgraph of
$\mathcal D_r$ (which has vertices ${\bf e}_i \pm {\bf e}_j$ for
$1\le i<j \le r$, where ${\bf e}_1$, \ldots, ${\bf e}_r$ is an
orthonormal set of vectors). From our knowledge of the structure
of cyclotomic signed graphs, we see that by relabelling and
changing signs of basis vectors (thereby inducing an equivalence),
we can suppose that $H$ contains ${\bf e}_1 + {\bf e}_2$, ${\bf
e}_2 + {\bf e}_3$, \ldots, ${\bf e}_{r-1} + {\bf e}_r$, and that
all other vertices of $H$ are of the form ${\bf e}_i - {\bf
e}_{i+1}$ (for some $i$ in the range $1\le i\le r-1$), or
${\bf e}_1 \pm {\bf e}_r$.

Now suppose that $w$ is a charged vertex in $G$. Since (i) $G$ is
connected, (ii) $Y_1$ is an excluded subgraph, and (iii) adjacent
charged vertices that have the same charge share all their
neighbours (Lemma \ref{L:localgeometry2}), we deduce that $w$ is
adjacent to one or more vertices in $H$. We treat first the case
where $w$ has charge $-1$. We have represented (a graph equivalent
to) $H$ by a set of Gram vectors, where adjacency of unequal
vertices is given by the dot product, and we can extend this to (a
graph equivalent to) $H\cup\{w\}$, where $w$ is represented by the
Gram vector
$ {\bf w} = \sum_{i=1}^{r+1} \lambda_i {\bf e}_i.$
If ${\bf w}$ is in the span of ${\bf e}_1$, \ldots, ${\bf e}_r$,
then we may set ${\bf e}_{r+1}={\bf 0}$; otherwise we need an
extra dimension for ${\bf w}$, and take ${\bf e}_{r+1}$ of length
1 and orthogonal to all of ${\bf e}_1$, \ldots, ${\bf e}_r$. Since
$w$ has charge $-1$, ${\bf w}$ has length $1$.

We consider two subcases. Case 1 (which we shall prove to be
impossible): $H$ contains one or both of ${\bf e}_1 \pm {\bf e}_r$,
so that $H$ contains a cycle of length $r$ containing no pair of
conjugate vertices. Case 2: $H$ contains neither of the vertices
${\bf e}_1 \pm {\bf e}_r$.

In Case 1, $H$ contains at least one cycle of length $r$ containing
no pair of conjugate vertices. Suppose that $w$ were adjacent to at
least two vertices on such a cycle, say $x$ and $y$ (and perhaps
others). Since subgraphs of $G$ equivalent to $Y_4$ have been
excluded, and $G$ cannot contain a subgraph equivalent to $X_1$, the
vertices $x$ and $y$ are not adjacent. By Lemma
\ref{L:localgeometry1}, every neighbour of $x$ is a neighbour of
$y$. But in a cycle of length at least 5 containing unadjacent
vertices $x$ and $y$ and containing no pair of conjugate vertices,
there will be a neighbour of $x$ that is not a neighbour of $y$.

Still in Case 1, suppose next that $w$ is adjacent to exactly one
vertex in some cycle of length $r$ containing no pair of conjugate
vertices. Then $G$ would contain a subgraph equivalent to $X_6$,
giving a contradiction.

To kill off Case 1, we now consider the remaining subcase where $w$
is adjacent to none of the vertices in the cycle ${\bf e}_1 + {\bf
e}_2$, ${\bf e}_2 + {\bf e}_3$, \ldots, ${\bf e}_{r-1} + {\bf e}_r$,
${\bf e}_1 \pm {\bf e}_r$. Then $H$ must contain at least one more
vertex, and after some relabelling and equivalence we can assume
that $w$ is adjacent to ${\bf e}_1 - {\bf e}_2$, with a positive
edge. Then
\[
\lambda_1 - \lambda_2 = 1\,,\quad \lambda_1 + \lambda_2 = \lambda_2
+ \lambda_3 = \lambda_3 + \lambda_4 = \lambda_4 \pm \lambda_5 = 0\,,
\]
where the `$\pm$' might be `$-$' if $r=5$. This gives
\[
\lambda_1 = 1/2\,,\quad \lambda_2 = -1/2\,,\quad \lambda_3 =
1/2\,,\quad \lambda_4 = -1/2\,,\quad \lambda_5 = \pm 1/2\,,
\]
and hence $|{\bf w}|>1$, giving a contradiction.

We now move to Case 2, where $H$ contains the path formed by the
vertices ${\bf e}_1 + {\bf e}_2$, ${\bf e}_2 + {\bf e}_3$, \ldots,
${\bf e}_{r-1} + {\bf e}_r$, and all other vertices in $H$ are of
the form ${\bf e}_i - {\bf e}_{i+1}$ (for some $i$ in the range
$1\le i\le r-1$).

If $w$ were adjacent to more than one vertex in our path, say $x$
and $y$, then as in Case 1 we would have $x$ and $y$ unadjacent,
implying that they share all their neighbours, giving a
contradiction.

If $w$ were not adjacent to any vertex in our path, then it would be
adjacent to some ${\bf e}_i - {\bf e}_{i+1}$, and from
\[
\lambda_i - \lambda_{i+1} = \pm 1\,,\quad \lambda_1 + \lambda_2 =
\lambda_2 + \lambda_3 = \ldots = \lambda_{r-1} + \lambda_r = 0\,,
\]
we would get at least five distinct $j$ such that $|\lambda_j|=1/2$,
contradicting $|{\bf w}|=1$.

We are reduced to the case where $w$ is adjacent to exactly one
vertex in our path. Since $X_6$ is excluded as a subgraph, this
neighbour of $w$ must be an endvertex of our path. Relabelling, we
can suppose that $w$ is attached to ${\bf e}_1 + {\bf e}_2$ by a
positive edge, but to none of ${\bf e}_2 + {\bf e}_3$, \ldots,
${\bf e}_{r-1} + {\bf e}_r$. If $H$ also contained ${\bf e}_1 -
{\bf e}_2$, then $w$ would necessarily be adjacent to it, or else
$G$ would contain a subgraph equivalent to $X_7$. Moreover, as
${\bf e}_2 + {\bf e}_3$ is joined to ${\bf e}_1 - {\bf e}_2$
by a negative edge, exclusion of subgraphs equivalent to $X_8$
implies that $w$ must then be connected to ${\bf e}_1 -
{\bf e}_2$ by a positive edge.

To sum up, if the minimal value of $r$ is at least 5, then we can
assume that $H$ contains the vertices ${\bf e}_1 + {\bf e}_2$, ${\bf
e}_2 + {\bf e}_3$, \ldots, ${\bf e}_{r-1} + {\bf e}_r$, and that all
other vertices are of the form ${\bf e}_i - {\bf e}_{i+1}$. Any
negatively charged vertex $w$ in $G$ is adjacent to one of ${\bf
e}_1 \pm {\bf e}_2$ or ${\bf e}_{r-1} \pm {\bf e}_r$. If both of
${\bf e}_1 \pm {\bf e}_2$  are in $H$ and $w$ is adjacent to one of
them, then it is adjacent to both; similarly for ${\bf e}_{r-1} \pm
{\bf e}_r$. The excluded graph $X_8$ constrains the signs of the
edges that connect $w$ to $H$. In short, $H\cup\{w\}$ is equivalent
to a subgraph of one of the $C_{2k}^{++}$ or $C_{2k}^{+-}$.

By equivalence, similar remarks hold for positively-charged vertices
in $G$.

If more than one charged vertex in $G$ is adjacent to the same
vertex in $H$, then the exclusion of subgraphs equivalent to
$Y_2$ and $Y_3$ implies that these charged vertices are
adjacent to each other; the exclusion of subgraph $Y_1$ implies
that they all have the same sign; Lemma \ref{L:localgeometry2}
implies that there are at most two such. We conclude that $G$ is
equivalent to a subgraph of one of the $C_{2k}^{++}$ or
$C_{2k}^{+-}$.

\subsubsection{$H$ is equivalent to a subgraph of $S_{16}$}
We shall show that $H$ is in fact equivalent to a subgraph of
$\mathcal D_r$ for some $r$, so that we are reduced to the
previous case.

Recalling previous notation, the vertices of $S_{16}$ are labelled
1, 2, 3, 4, 5, 6, 7, 8, $1234$, $1\bar{2}56$, $1\bar{3}\bar{5}7$,
$1\bar{4}\bar{6}\bar{7}$, $2\bar{3}58$, $2\bar{4}6\bar{8}$,
$3\bar{4}78$, $5\bar{6}7\bar{8}$ (a trivial relabelling of
\textbf{G4}). These are vectors in 8-dimensional real space,
with adjacency of unequal vectors given by the dot product.  Each
vector has length $\sqrt 2$.  Our restrictions on $G$ imply that
it has no triangles except perhaps involving two charged vertices
and one neutral vertex.

Note that $S_{16}$ is bipartite, with parts ${\mathcal
V}_1=\{1,2,3,4,5,6,7,8\}$, ${\mathcal V}_2=\{1234$, $1\bar{2}56$,
$1\bar{3}\bar{5}7$, $1\bar{4}\bar{6}\bar{7}$, $2\bar{3}58$,
$2\bar{4}6\bar{8}$, $3\bar{4}78$, $5\bar{6}7\bar{8}\}$.  There is
an equivalence of $S_{16}$ that interchanges these two parts,
induced by the orthogonal map with matrix
\[
\frac{1}{2}
\left(
\begin{array}{cccccccc}
1 & 1 & 1 & 1 & 0 & 0 & 0 & 0 \\
1 &-1 & 0 & 0 & 1 & 1 & 0 & 0 \\
1 & 0 &-1 & 0 &-1 & 0 & 1 & 0 \\
1 & 0 & 0 &-1 & 0 &-1 &-1 & 0 \\
0 & 1 &-1 & 0 & 1 & 0 & 0 & 1 \\
0 & 1 & 0 &-1 & 0 & 1 & 0 &-1 \\
0 & 0 & 1 &-1 & 0 & 0 & 1 & 1 \\
0 & 0 & 0 & 0 & 1 &-1 & 1 &-1
\end{array}
\right)
\]
with respect to ${\bf e}_1$, \ldots, ${\bf e}_8$.

Let $w$ be a charged vertex in $G$.  Arguing as before, $w$ is
adjacent to at least one vertex in $H$. First we treat the case
where $w$ is adjacent to at least two vertices in $H$, say $x$ and
$y$. Now $x$ and $y$ cannot be adjacent in $G$, or we would have a
forbidden triangle equivalent to $X_1$ or $Y_4$.  Then by
Lemma \ref{L:localgeometry1} the vertices $x$ and $y$ share all
their neighbours (and they must have at least one neighbour in $H$
or $H$ would not be connected). It follows that $x$ and $y$ are
either both in ${\mathcal V}_1$ or both in ${\mathcal V}_2$.
Working up to equivalence, and swapping ${\mathcal V}_1$ and
${\mathcal V}_2$ as above if necessary, we may suppose that $x$
and $y$ are both in ${\mathcal V}_1$. We may also suppose that $w$
is negatively charged, so that if we extend our set of Gram
vectors representing $H$ (some subset of the vectors/vertices in
${\mathcal V}_1\cup {\mathcal V}_2$) to a set of Gram vectors
representing $H\cup\{w\}$, the vector ${\bf w}$ representing $w$
will have length 1. We may write ${\bf w} = \sum_{i=1}^9 \lambda_i
{\bf e}_i\,$,
 where ${\bf e}_9$ (length $\sqrt{2}$, orthogonal to ${\bf e}_1$,
\ldots, ${\bf e}_8$) is included in case we need an extra
dimension to make room for ${\bf w}$. Since $|{\bf w}|=1$, we have
$\sum_{i=1}^9 \lambda_i^2 = 1/2\,.$

If $x$ and $y$ correspond to $i$ and $j$ in our labelling of the
vertices of $S_{16}$, then from
\begin{equation}\label{E-ij}
{\bf w}.{\bf e}_i=\pm 1,\quad {\bf w}.{\bf e}_j=\pm 1,
\end{equation}
we have $\lambda_i,\lambda_j\in\{1/2,-1/2\}$, and hence all other
$\lambda_k$ are zero.

There are now essentially two cases (up to equivalence):
$\{i,j\}=\{1,2\}$ and $\{i,j\}=\{1,8\}$.  Indeed there are
self-equivalences of $S_{16}$ induced by elements of the Klein
$4$-group acting on any of the six `missing' quartets
$\{1,2,7,8\}$, $\{1,3,6,8\}$, $\{1,4,5,8\}$, $\{2,3,6,7\}$,
$\{2,4,5,7\}$, $\{3,4,5,6\}$ (one then needs to apply appropriate
further transpositions of the form $(i\, \overline{i})$, and some
changes of signs of certain vertices, to map $S_{16}$ to itself).
We see that any vertex in ${\mathcal V}_1$ can be mapped to $1$ by
a self-equivalence of $S_{16}$, and that with $1$ fixed, any
vertex in ${\mathcal V}_1\backslash\{1,8\}$ can be mapped to $2$.

 In the case $\{i,j\}=\{1,2\}$,
 since $1$ and $2$ are not adjacent, Lemma \ref{L:localgeometry1}
implies that they have the same neighbours in $G$, and hence also
in $H$, whence $1\bar{3}\bar{5}7$, $1\bar{4}\bar{6}\bar{7}$,
$2\bar{3}58$, $2\bar{4}6\bar{8}$ $\not\in H$.  One of $1234$ and
$1\bar{2}56$ has dot product $\pm 1$ with $\bf w$, and hence must
be excluded from $H$ (or else together with $w$ and 1 (or 2) we would
have a forbidden triangle). Hence the vertices in $H$ are a
subset of $\mathcal W_1\cup \mathcal W_2$, where
\[
\mathcal W_1 = \{ 1,\ 2,\ 3,\ 4,\ 5,\ 6,\ 7,\ 8,\  3\bar{4}78,\
5\bar{6}7\bar{8}\},\quad  \mathcal W_2=\{1234\} \text{ or } \{1\bar{2}56\},
\]
depending on the signs of $\lambda_1$ and $\lambda_2$.  Then $H$
is readily seen to be equivalent to a subgraph of $\mathcal D_8$.

For the other essentially distinct case, $\{i,j\}=\{1,8\}$,
similar reasoning shows that $H$ is a subset of ${\mathcal
V}_1$, contradicting the connectedness of $H$.

We are left with the possibility that $w$ is adjacent to exactly
one vertex in $H$.  Let us temporarily call a signed charged graph
$K$ \emph{friendly} if it is cyclotomic,  contains exactly one
charged vertex $w$, the vertex $w$ is joined to exactly one
neutral vertex, and the neutral vertices in $K$ form a single
component.  In our current case, $H\cup\{w\}$ is friendly.  It
will be enough to show that any friendly graph with neutral
component equivalent to a subgraph of either $S_{16}$ or
$S_{14}$ is contained in a larger friendly graph (where the
neutral component of the larger friendly graph might or might not
be equivalent to a subgraph of either $S_{16}$ or $S_{14}$).
For then we can grow our friendly graph $H\cup\{w\}$ to a larger
friendly graph $H'\cup\{w\}$ with $H'$ not equivalent to a
subgraph of either $S_{16}$ or $S_{14}$. Then $H'$ must be
equivalent to a subgraph of some $\mathcal D_r$, and hence
the same is true for $H$.

A computer search checked that all friendly graphs with up to
$14$ neutral vertices are contained in larger friendly graphs.
As an indication of the work involved, some $377$ friendly
graphs with $15$ vertices ($14$ neutral vertices) were considered;
these would not all have been inequivalent, as it proved more
efficient to perform a fast but imperfect weeding out of
equivalent graphs, allowing some repeats through.  The search
could have been pushed further, but it was easier simply to check
 that there are no friendly graphs with $15$ or $16$
neutral vertices for which the neutral component is equivalent to
a subgraph of $S_{16}$.

\subsubsection{$H$ is equivalent to a subgraph of $S_{14}$}
The argument here is very similar to that for $S_{16}$, but in
fact slightly simpler, as $S_{14}$ has fewer vertices.
Analogously, we have $\mathcal V_1 = \{1, 2, 3, 4, 5, 6, 7\}$,
$\mathcal V_2 = \{1234, 1\bar{2}5\bar{6}, 1\bar{3}\bar{5}7,
1\bar{4}6\bar{7}, 2\bar{3}\bar{6}\bar{7}, 2\bar{4}57,
3\bar{4}\bar{5}\bar{6} \}$.  In the `unfriendly' case we find that
the vertices in $H$ are (after a suitable equivalence) a subset of
$\mathcal W_1 \cup \mathcal W_2$ where \[ \mathcal W_1 = \{1, 2,
3, 4, 5, 6 \},\qquad \mathcal W_2 = \{1234\}\ {\rm or}\
\{1\bar{2}5\bar{6}\}\,.
\] Then $H$ is equivalent to a subgraph of $\mathcal D_6$.

This completes the proof of Theorem \ref{T-LCSG}.


\section{Eigenvalues in the open interval $(-2,2)$}\label{S:open}

\subsection{Introduction to the next three sections} Sections \ref{S:open},\ref{S:charged} and \ref{S:nonneg} are devoted to results
for matrices and graphs under further restrictions. These follow
more or less straightforwardly from Theorems \ref{T-CSG} and
\ref{T-LCSG}. We consider first restricting to eigenvalues in the
open interval $(-2,2)$ (Section \ref{S:open}), deferring the
proofs to Section \ref{S:openproofs}. Then we consider charged
(unsigned) graphs, treating both the open and closed intervals
(Section \ref{S:charged}). Finally we treat symmetric matrices
that have non-negative integer entries (Section \ref{S:nonneg}).

\subsection{Cyclotomic signed graphs with all eigenvalues in $(-2,2)$.} Having classified all integer symmetric matrices having all their
eigenvalues in the interval $[-2,2]$, a natural question is what
happens if we restrict the eigenvalues to the open interval $(-2,2)$. From our
knowledge of the closed interval case, we can immediately restrict to
cyclotomic signed graphs and cyclotomic charged signed graphs, and need only consider
subgraphs of the maximal ones.

\begin{figure}[h]
\begin{center}
\leavevmode
\psfragscanon
\psfrag{U_1}{$U_1$}\psfrag{U_2}{$U_2$}\psfrag{U_3}{$U_3$}\psfrag{U_4}{$U_4$}
\psfrag{U_5}{$U_5$}\psfrag{U_6}{$U_6$}\psfrag{U_7}{$U_7$}\psfrag{U_8}{$U_8$}
\psfrag{U_9}{$U_9$}\psfrag{U_10}{$U_{10}$}\psfrag{U_11}{$U_{11}$}
\hbox{ \epsfxsize=4.3in
\epsffile{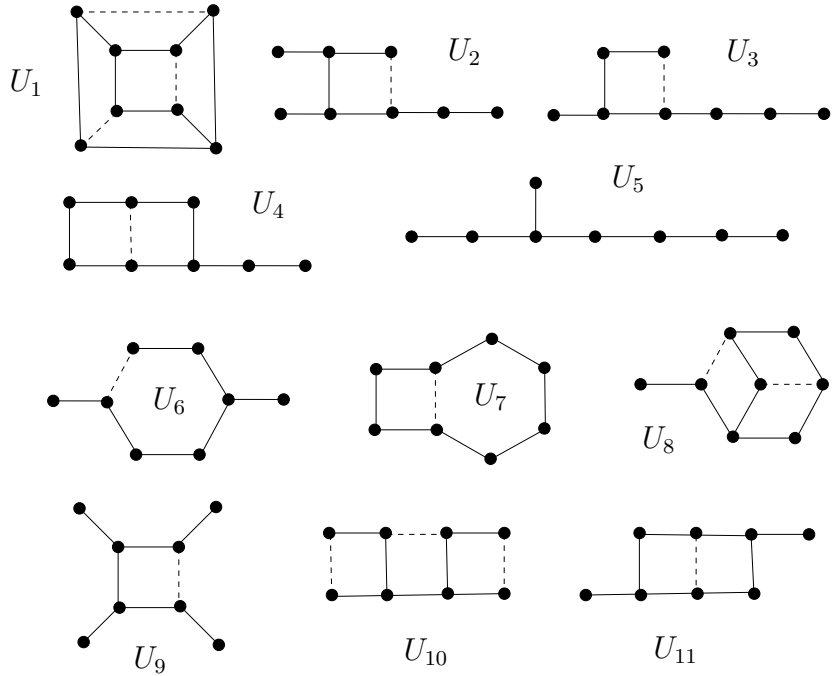} }
\end{center}
\caption{The sporadic connected cyclotomic signed graphs maximal
with respect to having all eigenvalues in $(-2,2)$.}  \label{F-U}
\end{figure}

\begin{figure}[h]
\begin{center}
\leavevmode
\psfragscanon
\hbox{ \epsfxsize=0.8in
\epsffile{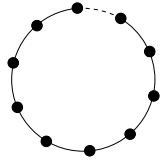} }
\end{center}
\caption{The $2k$-vertex connected cyclotomic signed graph $O_{2k}$, maximal
with respect to having all eigenvalues in $(-2,2)$,  shown here
for $k=5$.}  \label{F-O}
\end{figure}

\begin{figure}[h]
\begin{center}
\leavevmode
\psfragscanon
\psfrag{h}{$h$}\psfrag{k}{$k$}
\psfrag{b1}{$\overbrace{\phantom{XXXXXXXXXXXxxxx}}$}
\psfrag{b2}{$\overbrace{\phantom{XXXXXXXXXxxx}}$}
\hbox{ \epsfxsize=4.3in
\epsffile{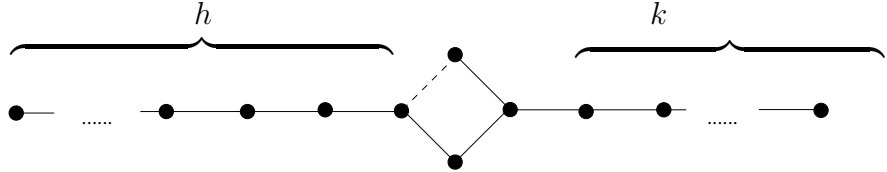} }
\end{center}
\caption{The doubly infinite family $Q_{hk}$ of connected
cyclotomic signed graphs
 having all eigenvalues in $(-2,2)$ but not contained in a maximal one.
}  \label{F-Q}
\end{figure}

\begin{theorem}[{``Uncharged, signed, $(-2,2)$''}]\label{T:open_signed}
Up to equivalence, the  connected signed  graphs maximal with
respect to having all their eigenvalues in $(-2,2)$ are the eleven
$8$-vertex sporadic examples $U_1,\dots,U_{11}$ shown in Figure
\ref{F-U}, and the infinite family $ O_{2k}$ of $2k$-cycles with
one edge of sign $-1$, for $2k\ge 8$, shown in Figure \ref{F-O}.

Further, every  connected cyclotomic signed graph having all its
eigenvalues in $(-2,2)$ is either contained in a maximal one, or
is a subgraph of one of the signed graphs $Q_{hk}$ of Figure
\ref{F-Q} for $h+k\ge 4$.
\end{theorem}

We note in passing that the graphs $U_i$ can all be obtained from
the cube $U_1$ by deleting certain edges. Not every choice of
edge-deletion produces a $U_i$, however. For instance no
edge-deleted subgraph of $U_1$ containing an induced subgraph
equivalent to $\tilde{D}_5$ can have all its eigenvalues in
$(-2,2)$.

\subsection{Cyclotomic charged signed graphs with all eigenvalues in $(-2,2)$.}

Next we have a corresponding result for charged signed graphs.

\begin{theorem}[{``Charged, signed, $(-2,2)$''}]\label{T:open_charged_signed}
Up to equivalence, the connected charged signed  graphs maximal
with respect to having all their eigenvalues in $(-2,2)$, and not
covered by the Theorem \ref{T:open_signed} above, are the  eight
$4$-vertex sporadic examples $V_1,V_2,\dots,V_8$ shown in Figure
\ref{F-V}, and the infinite family $P^{\pm}_n$ of $n$-vertex
charged paths of Figure \ref{F-P} for $n\ge 4$.

Further, every connected cyclotomic charged signed graph not
covered by the previous theorem is contained in such a maximal
one.
\end{theorem}

\subsection{Cyclotomic matrices with all eigenvalues in $(-2,2)$.}

We can combine the previous two theorems, translated into
matrix language, to obtain the following.

\begin{theorem}[{``Integer matrix, $(-2,2)$''}]\label{T-ISMopen} Every indecomposable  cyclotomic
matrix maximal with repect to having all its eigenvalues in the open
interval $(-2,2)$ is equivalent to the adjacency matrix of one of
the graphs $U_1, U_2,\dots,U_{11}$, $O_{2k}(2k\ge 8)$, $V_1,
V_2,\dots,V_8$, $P_n^\pm(n\ge 4)$ (given by Theorems
\ref{T:open_signed} and \ref{T:open_charged_signed}).

Further, every indecomposable  cyclotomic matrix  having all its
eigenvalues in $(-2,2)$ is either contained in a maximal one, or is
contained in the adjacency matrix of one of the signed graphs
$Q_{hk}$ of Figure \ref{F-Q} for $h+k\ge 4$.
\end{theorem}

\section{Maximal cyclotomic charged graphs}\label{S:charged}

\subsection{Cyclotomic charged unsigned graphs.} We now restrict our attention to cyclotomic charged graphs,
looking for those that are maximal with respect to having all
their eigenvalues in $[-2,2]$. For such a graph $G$ we need to
define $\overline{G}$ as the graph whose edges are the
 same as those of $G$, with the same signs, but with the charges on
 vertices being the opposite of those on $G$.  For example, graphs $V_1$ and $\overline{V_1}$ are shown in Figure \ref{F-V}.  It is clear that when $G$ is a
 tree, $\overline{G}$ is equivalent to $G$.

One difference for this kind of maximality is that it is not
 a property of equivalence classes of charged graphs: two of them may be equivalent with one of them maximal
 and the other not. For instance, one consequence of the next result is that
 for the maximal graphs $W_5$ and $W_6$ the graphs $\overline{W_5}$ and
 $\overline{W_6}$, being subgraphs of $W_7$, are not maximal.

 However, we have the following.

\begin{figure}[h]
\begin{center}
\leavevmode \psfragscanon
\psfrag{1}{$W_1$}\psfrag{2}{$W_2$}\psfrag{3}{$W_3$}\psfrag{4}{$W_4$}
\psfrag{5}{$W_5$}\psfrag{6}{$W_6$}\psfrag{7}{$W_7$}\psfrag{8}{$W_8$}
\psfrag{9}{$W_9$}\psfrag{10}{$W_{10}$}\psfrag{11}{$W_{11}$}
\psfrag{12}{$W_{12}$}\psfrag{13}{$W_{13}$}
\psfrag{1b}{$\overline{W_1}$}\psfrag{10b}{$\overline{W_{10}}$}\psfrag{11b}{$\overline{W_{11}}$}
\psfrag{12b}{$\overline{W_{12}}$} \hbox{ \epsfxsize=5.3in
\epsffile{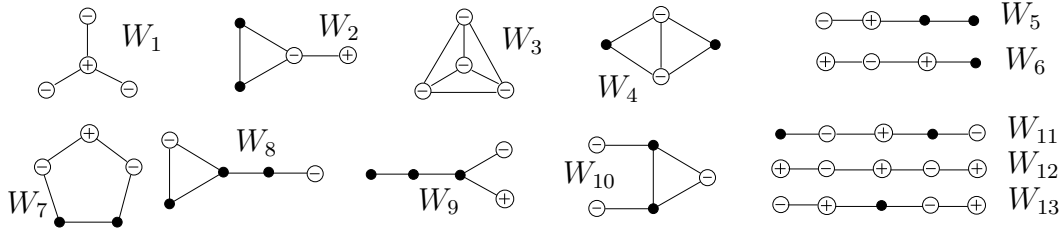} }
\end{center}
\caption{The sporadic maximal connected cyclotomic charged graphs
$W_1,\dots,W_{13}$. } \label{F-MCCS}
\end{figure}

\begin{figure}[h]
\begin{center}
\leavevmode \psfragscanon
\psfrag{F_5}{$F_5$}\psfrag{F_6}{$F_6$}\psfrag{F_n}{$F_n$}
\psfrag{G_5}{$G_5$}\psfrag{G_6}{$G_6$}\psfrag{G_n}{$G_n$}
\psfrag{H_3}{$H_3$}\psfrag{H_4}{$H_4$}\psfrag{H_n}{$H_n$}
\psfrag{I_3}{$I_3$}\psfrag{I_4}{$I_4$}\psfrag{I_n}{$I_n$}
\psfrag{J_2}{$J_2$}\psfrag{J_3}{$J_3$}\psfrag{J_n}{$J_n$} \hbox{
\epsfxsize=4.3in \epsffile{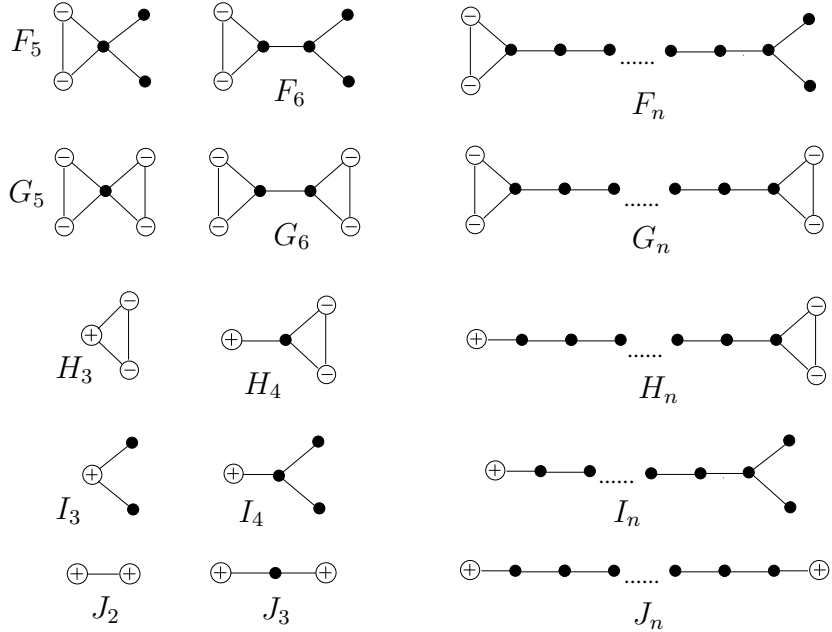} }
\end{center}
\caption{Five families of maximal connected cyclotomic charged
graphs. The smallest two members of each family are also shown. }
\label{F-MCCF}
\end{figure}

\begin{theorem}[{``Charged, unsigned, $[-2,2]$''}] \label{T:unsigned_charged}
The maximal connected cyclotomic charged  graphs
not covered by
Smith's result (i.e. not graphs),  are the sporadic examples
$W_1,\dots,W_{13}$ from Figure \ref{F-MCCS}, along with
$\overline{W_1}$,
$\overline{W_{11}}$,$\overline{W_{12}}$, and the seven families
$F_n(n\ge 5)$, $G_n(n\ge 5)$, $H_n(n\ge 3)$, $I_n(n\ge 3)$,
$J_n(n\ge 2)$ and $\overline{I_n}(n\ge 3)$, $\overline{J_n}(n\ge
2)$ from Figure \ref{F-MCCF}.

Further, every connected cyclotomic charged  graph  is contained
in such a maximal one.
\end{theorem}

The proof of this theorem is by inspection of the maximal
connected cyclotomic charged signed graphs of Theorem \ref{T-LCSG}
to find their maximal connected charged (unsigned) subgraphs.

\subsection{Cyclotomic charged unsigned graphs with all eigenvalues in $(-2,2)$.}
We next have a corresponding result for eigenvalues in the open
interval $(-2,2)$.

\begin{figure}[h]
\begin{center}
\leavevmode
\psfragscanon
\psfrag{bbb}{$\overline{V_1}$}\psfrag{V_1}{$V_1$}\psfrag{V_2}{$V_2$}
\psfrag{V_3}{$V_3$}\psfrag{V_4}{$V_4$}\psfrag{V_5}{$V_5$}\psfrag{V_6}{$V_6$}
\psfrag{V_7}{$V_7$}\psfrag{V_8}{$V_8$}
\hbox{ \epsfxsize=4.3in
\epsffile{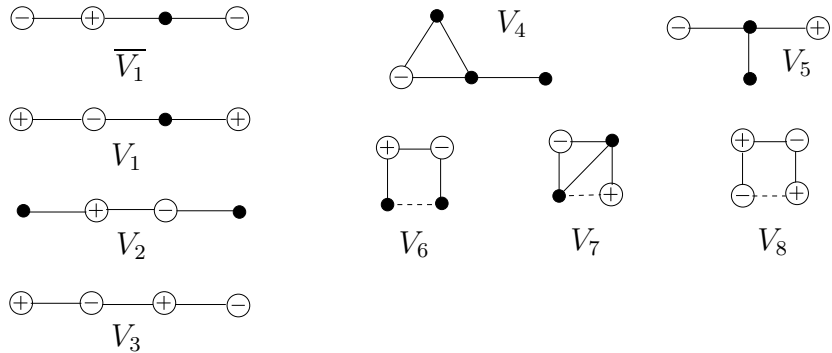} }
\end{center}
\caption{The sporadic connected cyclotomic charged signed graphs
maximal with respect to having all eigenvalues in $(-2,2)$.}
\label{F-V}
\end{figure}

\begin{figure}[h]
\begin{center}
\leavevmode
\psfragscanon
\psfrag{pn}{$P_n$}\psfrag{pn+}{$P_n^+$}\psfrag{pn-}{$P_n^-$}\psfrag{pn+-}{$P_n^{\pm}$}\psfrag{dn}{$D_n$}
\hbox{ \epsfxsize=2.2in
\epsffile{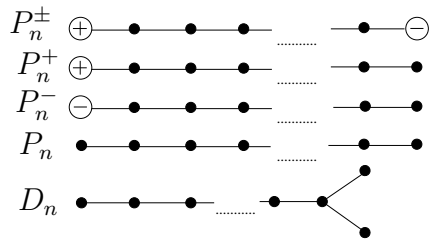} }
\end{center}
\caption{The $n$-vertex charged paths $P_n^{\pm}$ (for Theorems
\ref{T:open_charged_signed}, \ref{T-ISMopen} and
\ref{T:open_charged}), $P_n^+$ (Theorem
\ref{T:nonnegative_matrix_open}), $P_n^-$, $P_n$ (Section
\ref{S:polynomials}) and $D_n$ (Theorem \ref{T:nonnegative_matrix_open}).} \label{F-P}
\end{figure}

\begin{theorem}[{``Charged, unsigned, $(-2,2)$''}]\label{T:open_charged}
 The connected charged (unsigned)  graphs maximal with
respect to having all their eigenvalues in $(-2,2)$,   are the graph
$U_5$, the charged graphs
 $\overline{V_1}$, $V_1$, $V_2, \dots, V_5$, from Figure \ref{F-V}, and   $P^{\pm}_n$ of $n$-vertex charged paths
 of Figure \ref{F-P} for $n\ge 4$.

Further, every connected cyclotomic charged graph not covered by
Theorem \ref{T:open_signed} is contained in  one of the above
graphs.
\end{theorem}

Note that for $n\ge 8$ the graph $D_n$ (Figure \ref{F-P}) is a
subgraph of some $Q_{hk}$. So it has all its eigenvalues in
$(-2,2)$ and is covered by Theorem \ref{T:open_signed}. It is not
contained in any charged graph maximal with respect to having all
its eigenvalues in $(-2,2)$.

\section{maximal cyclotomic symmetric non-negative integer
matrices}\label{S:nonneg}

In this section we record our results for non-negative cyclotomic
matrices, i.e., those integer symmetric matrices that are
cyclotomic and have only non-negative entries.

\begin{theorem}[{``Non-negative integer matrix, $[-2,2]$''}] \label{T:nonnegative_matrix}
Up to conjugation by permutation matrices, the only maximal
indecomposable non-negative cyclotomic matrices
are the
matrices $(2)$ and $\left(\begin{matrix}0 & 2\\ 2 & 0
\end{matrix}\right)$, adjacency matrices of $\tilde{E}_6$,
$\tilde{E}_7$, $\tilde{E}_8$, $\tilde{A}_n (n\ge 2)$,
$\tilde{D}_n(n\ge 4)$ (Figure \ref{F-except}) along with the
two families  $I_n(n\ge 3)$ and $J_n(n\ge 2)$  (Figure
\ref{F-MCCF}).

Further, every  indecomposable non-negative cyclotomic matrix  is
contained in such a maximal one.
\end{theorem}

This result is readily deduced from Theorems
\ref{T-ISM},\ref{T:unsigned_charged} and Smith's results (Figure
\ref{F-except}).

\begin{theorem}[{``Non-negative integer matrix, $(-2,2)$''}] \label{T:nonnegative_matrix_open}
Up to conjugation by permutation matrices, the only indecomposable
non-negative cyclotomic matrix maximal with respect to having all
its eigenvalues in $(-2,2)$ is the adjacency matrix of $U_5$
(Figure \ref{F-U}).

Further, every  indecomposable non-negative cyclotomic matrix  is
either contained in the adjacency matrix of $U_5$ or in the
adjacency matrix of either $P_n^+$ or $D_n$ (Figure \ref{F-P}) for
some $n$.
\end{theorem}

\section{Proofs of Theorems \ref{T:open_signed} and
\ref{T:open_charged_signed}}\label{S:openproofs}
 To prove
Theorem \ref{T:open_signed}, we first we show that the two
infinite families $O_{2k}$, $Q_{hk}$ have their eigenvalues in the
open interval.

Suitable sets of Gram vectors, for the two cases, are:

For $O_{2k}$, the columns of the $(2k)\times(2k)$ matrix $(c_{ij})$, where
$$
c_{ij}=\begin{cases}
 1 \text{ if } i=j \text{ or } i=j+1\\
 -1 \text{ if } (i,j)=(1,2k)\\
 0 \text{ otherwise }.\\
\end{cases}
$$

For $Q_{hk}$, the columns of the $(h+k+4)\times(h+k+4)$ matrix $(q_{ij})$, where
$$
q_{ij}=\begin{cases}
 1 \text{ if } i=j\ (j=1,\dots,h+k+4) \text{ or } i=j+1\ (j=1,2,3) \text{ or } (i,j)=(2,k+5)\\
 \quad \text{ or } i=j-1\ (j=5,\dots,k+4,k+6,\dots,h+k+4)\\
 -1 \text{ if } (i,j)=(1,4)\\
 0 \text{ otherwise. }\\
\end{cases}
$$

Note that both of these sets of columns are easily seen to be
linearly independent. Hence, in each case, for the adjacency matrix
$A$ of these signed graphs,  $A+2I$ is non-singular, so $-2$ is not
an eigenvalue. Since these families comprise bipartite graphs (in
the extended sense), $2$ is not an eigenvalue. The families are
cyclotomic, being subgraphs of $T_n$ for some $n$, so we are done.

Now we find the remaining graphs.

For subgraphs of the sporadic graphs $S_{14}$ and $S_{16}$, we
know that these are subgraphs of $\mathcal E_8$, and so can be
embedded in $\mathbb R^8$, with $A+2I$ nonsingular. Hence such a
subgraph can have at most $8$ vertices. These can be found by
exhaustive search; the maximal ones are $U_1,\dots,U_{11}$ and
$O_8$.

There remain the subgraphs of the infinite families.

We observe that:
\begin{itemize}
\item an hour-glass \hbox{ \epsfxsize=0.3in
\epsffile{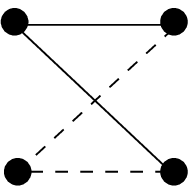} },
 equivalent to an unsigned square,  has $2$ as an eigenvalue;

\item the classical $\tilde D_n$ graphs (see Figure \ref{F-except})
have $2$ as an eigenvalue.
\end{itemize}

Hence the subgraph can contain at most one pair of conjugate
vertices. So it is either a path, a cycle,  some $Q_{hk}$ or
$Q'_{k}$, defined to be $Q_{1k}$ with its two leaves identified. A
path is a subgraph of some $Q_{hk}$, while a cycle must be
equivalent to some $O_{2k}$, for otherwise it is equivalent to a
cycle with all positive edges for which $2$ is an eigenvalue. For
$Q'_{k}$, we can delete one of its pair of conjugate vertices
to obtain a graph equivalent to a cycle with all positive edges.

This completes the proof of Theorem \ref{T:open_signed}.

The proof of Theorem \ref{T:open_charged_signed} is similar. We
can assume that the charged  graphs we seek do indeed have at
least one charged vertex. The relevant subgraphs of $S_7$, $S_8$
and  $S_8'$ are found by exhaustive search. For the subgraphs of
$C_{2k}^{++}$ and  $C_{2k}^{+-}$, we see by the same argument as
above that the neutral component can contain at most one pair
of conjugate vertices. Hence the neutral component is a path, or some $Q_{hk}$.

Two adjacent charges of the same sign have one of $\pm 2$ as an
eigenvalue, so each charged component has exactly one charge.
Putting charges of the same sign at each end of a path would give
one of $\pm 2$ as an eigenvalue, as one can see by writing down an
obvious eigenvector.

Putting a charge on either end of some $Q_{hk}$ gives one of $\pm
2$ as an eigenvalue. To see this  it suffices to consider adding a
negative charge to one end, with corresponding column vector
$(0,\dots,0,1)^T$ to add to $(q_{ij})$,  and adding a row of
zeroes to $(q_{ij})$ to make it square, giving a singular matrix,
and hence $-2$ as an eigenvalue.
 This leaves
$P_n^\pm$ (and its subgraphs) as the only possibilities.

For $P_n^\pm$, the columns of the $n\times n$ matrix $(p_{ij})$, where
$$
p_{ij}=\begin{cases}
  \sqrt{2} \text{ if } (i,j)=(1,1)\\
  1 \text{ if } i=j   \text{ or } i=j+1\quad (i\ge 2)\\
 0 \text{ otherwise\,,}\\
\end{cases}
$$
are easily seen to be linearly independent.
Again, for the adjacency matrix $A$ of this bipartite charged signed graph,  $A+2I$ is non-singular, so $-2$ is not an eigenvalue, and hence neither is $2$.

This completes the proof of Theorem \ref{T:open_charged_signed}.
Theorem \ref{T:open_charged} then follows easily.

\section{The cyclotomic polynomials of charged signed
graphs}\label{S:polynomials}
 Table \ref{Ta-1} gives the reciprocal
polynomials of the maximal connected cyclotomic charged signed
graphs that appear in our results. All are maximal in the sense
explained where they appear, apart from the $Q_{hk}$ which, as we
have seen, do not belong to any connected cyclotomic charged
signed graph maximal with respect to having all eigenvalues in
$(-2,2)$. Note, however, that the polynomials associated to
$C^{++}_{2k}$ and $S_7$  will need changes of variable $x\mapsto
-x$,  $z\mapsto -z$ when going from one equivalent, but not
strongly equivalent, graph to another.

\begin{table}[ht]
\begin{center}
\bigskip\begin{tabular} { c | c | c}

Charged signed graph  & Characteristic polynomial & Associated cyclotomic polynomial\\
     \hline
     $T_{2k}$ & $(x+2)^{k}(x-2)^k$ & $(z^2-1)^{2k}\quad (k\ge 3)$\\
     $S_{14}$ & $(x+2)^{7}(x-2)^7$ & $(z^2-1)^{14}$\\
 $S_{16}$ & $(x+2)^{8}(x-2)^8$ & $(z^2-1)^{16}$\\
$C^{++}_{2k}$ & $(x+2)^{k-1}(x-2)^{k+1}$ & $(z-1)^{2k+2}(z+1)^{2k-2}\quad (k\ge 2)$\\
$C^{+-}_{2k}$ & $(x+2)^{k}(x-2)^k$ & $(z^2-1)^{2k}\quad (k\ge 2)$\\
$S_7$ & $(x+2)^{3}(x-2)^{4}$ & $(z+1)^6(z-1)^8$\\
$S_8$, $S'_8$ & $(x+2)^{4}(x-2)^{4}$ & $(z^2-1)^8$\\
   \end{tabular}
    \vspace*{1ex}
    \caption{The characteristic and cyclotomic polynomials of maximal cyclotomic charged signed graphs.} \label{Ta-1}
\end{center}\end{table}

 Table \ref{Ta-2} gives the reciprocal polynomials of the cyclotomic signed graphs of Theorems
 \ref{T:open_signed} and \ref{T:open_charged_signed}, shown in Figures \ref{F-U},
\ref{F-Q}, \ref{F-V}, \ref{F-O} and \ref{F-P}.   In the table,
$\Phi_n$ denotes the $n$th cyclotomic polynomial.

\begin{table}[ht]
\begin{center}
\bigskip\begin{tabular} { c | c }

Charged signed graph  &  Associated cyclotomic polynomial\\
     \hline
$O_{2k}$ & $(z^{2k}+1)^2$\\
$Q_{hk}$ & $(z^{2h+4}+1)(z^{2k+4}+1)\quad (h+k\ge 4)$\\
   $U_1$  & $\Phi_{6}^4(z^2)$\\
        $U_2$   & $\Phi_{20}(z^2)$\\
    $U_3$ & $\Phi_{24}(z^2)$\\
    $U_4$ & $\Phi_{6}(z^2)\Phi_{18}(z^2)$\\
    $U_5$ & $\Phi_{30}(z^2)$\\
     $U_6$, $U_9$   & $\Phi_{12}^2(z^2)$\\
    $U_7$, $U_{11}$ & $\Phi_{15}(z^2)$\\
 $U_8$ & $\Phi_{12}(z^2)\Phi_{6}^2(z^2)$\\
    $U_{10}$ & $\Phi_{10}^2(z^2)$\\
      $V_1$        & $\Phi_{15}(z)$\\
    $\overline{V_1}$,  $V_4$  & $\Phi_{30}(z)$\\
       $V_3$, $V_6$   & $\Phi_{20}(z)$\\
   $V_2$, $V_5$  & $\Phi_{24}(z)$\\
    $V_7$, $V_8$   & $\Phi_{12}^2(z)$\\
       $P_n$ & $(z^{2n+2}-1)/(z^2-1)$\\
$P^-_n$ & $(z^{2n+1}-1)/(z-1)$\\
$P^\pm_n$ & $z^{2n}+1\quad (n\ge 4)$\\
   \end{tabular}
    \vspace*{1ex}
    \caption{The cyclotomic polynomials of some charged signed graphs
    having all their eigenvalues in $(-2,2)$.} \label{Ta-2}
\end{center}\end{table}

For a single sporadic graph with adjacency matrix $A$, the
reciprocal polynomial $z^n\chi_A(z+1/z)$ can be easily calculated.
For the infinite families, more work is required. Here, for
convenience, we use the same notation for a graph and its
associated cyclotomic polynomial.


For computing formulae for families of associated cyclotomic
polynomials, a standard tool will be to use induction on the
determinant $\det((z+1/z)I-A)$, where $A$ is the adjacency matrix
of the graph under consideration. In this way it is easy first to
compute the the $n$-vertex (unsigned) path $P_n$, giving
$P_n=(z^{2n+2}-1)/(z^2-1)$, as in the table (see also \cite{MS}).
Then expansion by the first row of the determinant gives
$P_n^-=(z^{2n+1}-1)/(z-1)$. Also  $P_n^\pm$ is readily calculated,
again expanding in the same way.


For $O_{2k}$, determinant expansion firstly along the top row, and
then down the left rows of the resulting determinants gives
$O_{2k}=(z^2+1)P_{2k-1}-2z^2P_{2k-2}+2z^{2k}$, and hence the
result.

For $Q_{hk}$, the formulae for $Q_{1k}$ and $Q_{2k}$ can be proved
by induction, using the determinant, in a similar way to that for
$P_n^-$. These can then
be used as the base cases for an inductive proof of the $Q_{hk}$
formula.

For $T_{2k}$, label its top vertices $1,3,5,\dots,2k-1$ and the
bottom vertices $2,4,6,\dots,2k$, with $2$ the conjugate vertex to
vertex $1$. Then $(-1,1,1,1,0,\dots,0)$ is an eigenvector of
$T_{2k}$ with eigenvalue $-2$, and $(1,-1,1,1,0,\dots,0)$ is an
eigenvector of $T_{2k}$ with eigenvalue $2$, both associated to
the hourglass $[1,2,3,4]$. From the symmetry of $T_{2k}$ that acts
by $i\mapsto i+2 \mod 2k$ on its vertices, we get two
eigenvectors, with eigenvalues $-2$ and $2$ for each of the
hourglasses $[3,4,5,6], [5,6,7,8],\dots,[2k-1,2k,1,2]$.  These
eigenvectors are independent, so that $T_{2k}$ has characteristic
polynomial $(x+2)^k(x-2)^k$, which,
 on putting $x=z+1/z$, gives the result.

For $C_{2k}^{++}$, label the vertices as for $T_{2k}$.
For the hourglasses $[3,4,5,6]$, $\dots,[2k-5,2k-4,2k-3,2k-2]$
(those without charged vertices), we  get the same eigenvectors as
for $T_{2k}$, with the same eigenvalues. The hourglasses
$[1,2,3,4]$ and $[2k-3,2k-2,2k-1,2k]$ give the same eigenvectors
as for $T_{2k}$  with eigenvalue $-2$. For the hourglass
$[1,2,3,4]$, however, we also get two independent eigenvectors
$(1,1,0,\dots,0)$ and $(2,0,1,1,0,\dots,0)$ with eigenvalue $2$,
and  from the hourglass $[2k-3,2k-2,2k-1,2k]$ we get two more
independent eigenvectors $(0,\dots,0,-1,1)$ and
$(0,\dots,0,1,-1,0,2)$ with eigenvalue $2$. Thus $C_{2k}^{++}$ has
characteristic polynomial $(x+2)^{k-1}(x-2)^{k+1}$, giving the
result.

 For $C_{2k}^{+-}$, note that this is bipartite in the extended
sense, so that the eigenvalues $2$ and $-2$ have equal
multiplicities.

\section{Final remarks}
\subsection{Finite reflection groups}

Given the root system $\Phi$ of a finite reflection
group, one classically looks for a subset
$\Delta$ that is a \emph{simple system},
namely one that is a basis for the $\mathbb R$-span
of $\Phi$ and such that every element of $\Phi$
is a linear combination of elements of $\Delta$
with all coefficients weakly of the same sign.
The Coxeter graph of a simple system is determined
by the reflection group, and provides a means
of classifying finite reflection groups.

If we have a signed graph with all eigenvalues in
$(-2,2)$ then (as we have seen) its vertices
can be associated with a linearly independent set
$\Delta'$ of vectors, and the reflection group
generated by the hyperplanes orthogonal to those
vectors is a finite reflection group.
The closure of $\Delta'$ under this reflection
group is a root system $\Phi$:
in the language of \cite{CvL} we are
taking the star closure of the lines spanned by
the elements of $\Delta'$.

Our set $\Delta'$ will not generally be a simple
system for $\Phi$, but it will be a basis
for the $\mathbb R$-span of $\Phi$.
The unsigned version of our graph (making
all edges positive) is the Coxeter graph of
$\Delta'$.

The neutral signed graphs of Theorem \ref{T:open_signed} therefore
provide a classification of all Coxeter graphs coming from bases
for the $\mathbb R$-span of root systems contained in either
$\mathcal D_n$ $(n\ge 4)$ or $\mathcal E_8$. For example,
one can generate $\mathcal E_8$ using eight reflections whose
Coxeter graph is the  cube $U_1$ of Figure \ref{F-U}.

For other connections between signed graphs and Coxeter graphs and
roots systems see \cite{CST} and \cite {Z1}.

\subsection{The graph $\mathbf {S_{14}}$}\label{S-Chapman}
Robin Chapman has pointed out that, up to equivalence,  the signed
graph $S_{14}$ of Figure \ref{F-CSG_14} can be defined as follows:
label the vertices $0,1,\dots,6,0',1',\dots,6'$ and, working
modulo $7$, for each $i$ join $i$  to each of $i'$, $(i+1)'$ and
$(i+3)'$ by  positive edges, and join $i$ to $(i-1)'$ by a
negative edge. The representation of $S_{14}$ in the figure is
based on this observation.

\subsection{Chebyshev polynomials and cyclotomic matrices}
Let $\mathcal T_n(x)$ denote the $n$th Chebyshev polynomial of the first kind, defined on the interval $[-2,2]$. So $\mathcal T_n(x)$ has integer coefficients and satisfies
\begin{equation}\label{E-Ch}
\mathcal T_n\left(z+\frac{1}{z}\right)=z^n+\frac{1}{z^n}.
\end{equation}

Then for any cyclotomic matrix $A$, the matrix $\mathcal T_n(A)$
is again cyclotomic. This follows from diagonalizing $A$ and using
(\ref{E-Ch}).

\subsection{Acknowledgments}
The authors are grateful for the hospitality provided by the
University of Bristol during the time that this paper was written.
We thank John Byatt-Smith for producing Figure \ref{F-torus}.

\end{document}